\documentclass[twocolumn]{autart1}    

\usepackage{amssymb,amsmath}
\usepackage{enumerate,array,slashbox}
\usepackage{graphicx}          

\newtheorem{assumption}{Assumption}

\newtheorem{definition}{Definition}

\begin{document}

\begin{frontmatter}

\title{Newton based Stochastic Optimization using 
\\$q$-Gaussian Smoothed Functional Algorithms} 

\author{Debarghya Ghoshdastidar}\ead{debarghya.g@csa.iisc.ernet.in},   
\author{Ambedkar Dukkipati}\ead{ad@csa.iisc.ernet.in},              
\author{Shalabh Bhatnagar}\ead{shalabh@csa.iisc.ernet.in}  

\address{Department of Computer Science \& Automation,\\
Indian Institute of Science, Bangalore.}
          
\begin{keyword}                           
Smoothed functional algorithms, $q$-Gaussian perturbations,
Hessian estimate, stochastic optimization, two-timescale algorithms.
\end{keyword}                             

\begin{abstract}                          
We present the first $q$-Gaussian smoothed functional (SF) estimator of the Hessian and 
the first Newton-based stochastic optimization algorithm that estimates both the 
Hessian and the gradient of the objective function using $q$-Gaussian perturbations.
Our algorithm requires only two system simulations (regardless of the parameter dimension)
and estimates both the gradient and the Hessian at each update epoch using these.
We also present a proof of convergence of the proposed algorithm.
In a related recent work (Ghoshdastidar et al., 2013), we presented gradient SF
algorithms based on the $q$-Gaussian perturbations. Our work extends prior work on 
smoothed functional algorithms by generalizing the class of perturbation distributions
as most distributions reported in the literature for which SF algorithms are known 
to work and turn out to be special cases of the $q$-Gaussian distribution.
Besides studying the convergence properties of our algorithm analytically,
we also show the results of several numerical simulations on a model of a 
queuing network, that illustrate the significance of the proposed method. 
In particular, we observe that our algorithm performs better in most cases,
over a wide range of $q$-values, in comparison to Newton SF algorithms 
with the Gaussian (Bhatnagar, 2007) and Cauchy perturbations, as 
well as the gradient $q$-Gaussian SF algorithms (Ghoshdastidar et al., 2013).

%
\end{abstract}

\end{frontmatter}

\section{Introduction}

Stochastic techniques for optimization have gained immense popularity
over the last couple of decades. 
Stochastic alternatives have been developed for a variety of
classic optimization problems, 
such as maximum likelihood estimation~\cite{Spall_1987_conf_ACC}, 
expectation maximization~\cite{Delyon_1999_jour_AnnStat}, 
least squares estimation~\cite{Ljung_2001_jour_ACSP}, 
discrete parameter optimization~\cite{Mishra_2007_conf_CDC} and
control of discrete-event systems~\cite{Bhatnagar_1998_jour_ProbEnggInfoSc}
among others.
On the other hand,
a wide class of problems related to automated control and sequential decision making 
are often posed as Markov Decision Process (MDP) models~\cite{Puterman_1994_book_Wiley}.
A classic example is encountered in control of a stochastic process, 
where one needs to optimize the performance of a system by appropriately tuning some parameter.
This scenario is quite common in reinforcement learning problems~\cite{Bertsekas_1996_book_Athena}.
Simulation based schemes~\cite{Bhatnagar_2013_book_Springer}
are popularly used for solving MDPs since these algorithms do not require prior knowledge of
the system dynamics; rather, the transitions of the system are simulated to obtain
estimates of the cost function to be minimized.
 
One of the earliest ideas of stochastic optimization is due to 
Kiefer and Wolfowitz~\cite{Kiefer_1952_jour_AnnMathStat}, where the zeros
of the gradient of the objective function is determined via
Robbins-Monro root-finding~\cite{Robbins_1951_jour_AnnMathStat}.
This approach uses a finite difference gradient estimate,
and hence, is termed as finite difference stochastic approximation (FDSA).
It proves quite useful in optimization of 
stochastic functions, commonly encountered in stochastic control problems.
However, it requires $2N$ parallel simulations 
to estimate the gradient 
at each iteration, $N$ being the dimension of the optimizer.
More efficient techniques have been proposed in the literature, 
which perform gradient estimation using only two parallel simulations.
These techniques include 
simultaneous perturbations stochastic approximation (SPSA)~\cite{Spall_1992_jour_AutoControlTrans},
random direction stochastic approximation (RDSA)~\cite{Kushner_1978_book_Springer}
and the smoothed functional (SF) method~\cite{Styblinski_1990_jour_NeuNet}.
More computationally efficient one-simulation variants of these methods also 
exist~\cite{Styblinski_1986_jour_TransCADICS,Spall_1997_jour_Automatica}
but their performance is relatively poor when compared with their two-simulation counterparts.

All the above approaches employ an approximate \emph{steepest} descent method for optimization.
It is well known that second order techniques such as Newton's method, 
are faster and provide greater accuracy when compared with steepest descent methods.
In the context of stochastic optimization,
Newton based approaches have been proposed for FDSA~\cite{Ruppert_1985_jour_AnnStat},
SPSA~\cite{Spall_2000_jour_AutoCtrl,Bhatnagar_2005_jour_TOMACS}, and SF~\cite{Bhatnagar_2007_jour_TOMACS}, respectively,
in the literature. Though these methods suffer from increased 
computational effort due to Hessian estimation, 
projection (to the set of positive definite and symmetric matrices)
and inversion at each update,
they have been observed to perform significantly better
than their steepest descent counterparts~\cite{Spall_2000_jour_AutoCtrl,Bhatnagar_2007_jour_TOMACS}.
Approximate methods for projecting the Hessian to the set of positive definite matrices and its inversion have also been 
studied~\cite{Zhu_2002_jour_AdapCtrlSP,Bhatnagar_2005_jour_TOMACS}. 
An efficient approximation is the Jacobi variant, where only the diagonal
terms of the Hessian matrix are updated, while the off-diagonal elements are simply set to zero.
This simplifies both the Hessian projection 
and inversion procedures.

In this paper, we focus on the SF algorithms for optimization.
The idea of smoothing dates back to Katkovnik and Kulchitsky~\cite{Katkovnik_1972_jour_Automation}
and Rubinstein~\cite{Rubinstein_1981_book_Wiley},
where it was shown that a smoother variant of the objective function 
can be estimated as an expectation (sample average is used for practical purposes)
of perturbed observations of the objective function.
This method was employed by Styblinski and Opalski~\cite{Styblinski_1986_jour_TransCADICS}
for parameter optimization in the manufacturing process of integrated circuits.
Subsequently, a two-sided version of SF~\cite{Rubinstein_1981_book_Wiley}
was employed for gradient based optimization in~\cite{Styblinski_1990_jour_NeuNet}.
Second-order SF schemes were proposed by Bhatnagar~\cite{Bhatnagar_2007_jour_TOMACS}.
Till this point, SF methods considered smoothing using either Gaussian or Cauchy
distributions~\cite{Styblinski_1990_jour_NeuNet}, while it was known that
uniform and symmetric Beta distributions were possible
candidates~\cite{Rubinstein_1981_book_Wiley,Kreimer_1992_jour_AnnOR}. 
Moreover, the random search algorithms were also observed to be a special case of SF~\cite{Rubinstein_1981_book_Wiley}.
Recent studies by the authors~\cite{Ghoshdastidar_2012_conf_ISIT,Ghoshdastidar_2013_arxiv}
revealed that this set of smoothing kernels can be further
extended to the class of $q$-Gaussian distributions~\cite{Prato_1999_jour_PhyRevE},
extensively studied in the field of nonextensive statistical mechanics~\cite{Tsallis_1998_jour_PhysicaA}.

The main focus of the work in~\cite{Ghoshdastidar_2013_arxiv}
was to exploit two key properties of $q$-Gaussians. This class of 
distributions generalize the Gaussian distribution via a Tsallis generalization
of the Shannon entropy functional~\cite{Tsallis_1998_jour_PhysicaA}. 
Hence, they retain some of the \emph{nice} characteristics
of Gaussian, which include its smoothing properties~\cite{Rubinstein_1981_book_Wiley}.
On the other hand, the $q$-Gaussians exhibit a power-law nature for some $q$-values,
which can be exploited to achieve greater exploration in the SF method.
In~\cite{Ghoshdastidar_2013_arxiv}, the authors proposed two gradient descent optimization
algorithms using one-simulation and two-simulation $q$-Gaussian SF.
The current paper extends the work in~\cite{Ghoshdastidar_2013_arxiv} to Newton based search techniques
and derives, for the first time, a $q$-Gaussian Hessian estimator.
Due to the observation that two-simulation methods
consistently perform better than one-simulation 
methods~\cite{Bhatnagar_2005_jour_TOMACS,Bhatnagar_2007_jour_TOMACS,Ghoshdastidar_2013_arxiv},
we focus only on the two-simulation version of Newton based optimization
using $q$-Gaussian SF even though the one-simulation version is easier to derive 
than the two-simulation estimator that we present. 
We refer to our algorithm as N$q$-SF2
to indicate that it is a Newton based smoothed functional algorithm 
that uses $q$-Gaussian perturbations and requires two simulations.
This terminology is also consistent with that used in~\cite{Bhatnagar_2007_jour_TOMACS,Ghoshdastidar_2013_arxiv}.
Our approach requires estimation of the Hessian of the two-sided $q$-Gaussian SF.
We derive the Hessian estimator and subsequently present a 
theoretical analysis, which shows that the N$q$-SF2
algorithm converges to the neighborhood of a local optimum.
For various values of $q$, the $q$-Gaussian distribution encompasses the Gaussian,
Cauchy, symmetric Beta and uniform distributions as special cases.
Our analysis shows that the algorithm converges for a range of values of 
$q$ that includes the aforementioned distributions barring the uniform distribution.
Thus our work significantly enhances the class of perturbation distributions 
for smoothed functional algorithms, We note here that the Hessian estimator in~\cite{Bhatnagar_2007_jour_TOMACS}
is only for the case of Gaussian perturbations.
Simulations on a two-node queuing network show significant performance 
improvements of the N$q$-SF2 algorithm over the Newton algorithm in~\cite{Bhatnagar_2007_jour_TOMACS}
and also the gradient $q$-SF algorithms in~\cite{Ghoshdastidar_2013_arxiv}.

The rest of the paper is organized as follows:
Section~\ref{sec_hessian} briefly reviews the basic idea of smoothed functional (SF)
algorithms using the multivariate $q$-Gaussian distribution.
It also presents the $q$-Gaussian smoothed Hessian estimates that we use for our
Newton based search techniques.
Section~\ref{sec_problem} presents the proposed two-simulation method. 
This section describes the problem setting, and
gives the proposed algorithm.
Section~\ref{sec_convergence} provides a theoretical study of the 
convergence analysis of the algorithm,
while results of numerical experiments are shown in Section~\ref{sec_simulation}.
Concluding remarks are provided in Section~\ref{sec_conclusion}.
Finally, Appendix~\ref{app_convergence} provides some of the detailed
proofs pertaining to the convergence analysis.


\section{Hessian estimation using $q$-Gaussian SF}
\label{sec_hessian}
The idea of smoothed functionals (SF) was first proposed in~\cite{Katkovnik_1972_jour_Automation}.
Consider the optimization problem
\begin{align}
 &\min_{\theta\in C \subset \mathbb{R}^N} J(\theta),
\end{align}
where $C$ is a compact and convex subset of $\mathbb{R}^N$,
and $J:C\to \mathbb{R}$ is a real-valued function,
which either does not have an analytic expression,
or has an expression that is not known.
To achieve a ``good'' optimal solution in such cases, 
it is often useful to minimize a smoothed variant of the objective function,
called the smoothed functional, defined as
\begin{align}
 S_{\beta} [J(\theta)] 
 = \int_{\mathbb{R}^N} G_\beta(\eta) J(\theta-\eta) \mathrm{d}\eta \;.
\end{align}
Katkovnik and Kulchitsky~\cite{Katkovnik_1972_jour_Automation} 
considered $G_\beta$ to be the $N$-dimensional Gaussian distribution with zero mean
and covariance matrix $\beta^2 I_{N\times N}$.
An alternative (two-sided) definition of the smoothed functional
was given in~\cite{Rubinstein_1981_book_Wiley} as
\begin{align}
 S_{\beta} [J(\theta)] 
  &= \frac{1}{2} \int_{\mathbb{R}^N}  G_\beta(\eta) \big(J(\theta+\eta) 
  + J(\theta-\eta) \big) \mathrm{d}\eta \;.
  \label{eq_2SF_general}
\end{align}
It was shown that reasonably good solutions can be obtained
by minimizing $S_\beta[J(\cdot)]$ using standard optimization techniques 
with properly tuned smoothing parameter $\beta$.
Subsequently, Rubinstein~\cite{Rubinstein_1981_book_Wiley} showed that
Gaussian distribution is not the only such function,
and provided necessary conditions that need to be satisfied by
a distribution to provide appropriate smoothing. 
The uniform and Cauchy distributions were also found to satisfy
the properties required of smoothing kernels.
A symmetric version of the Beta distribution was shown to satisfy 
similar conditions~\cite{Kreimer_1992_jour_AnnOR}, 
and this smoothing was shown to have connections with polynomial approximations.

\subsection{$q$-Gaussian smoothed functionals} 

Recently, Ghoshdastidar et al.~\cite{Ghoshdastidar_2013_arxiv}
proposed a class of smoothing kernels based on 
the $N$-dimensional $q$-Gaussian distributions~\cite{Tsallis_2007_arxiv,Vignat_2007_jour_PhyA} defined as
\begin{align}
\label{Gq:formula_q}
G_{q}(&x|q,\mu_q,\Sigma_q) = {\frac{1}{K_{q,N} |\Sigma_q|^{1/2}}} \times
\\&
\left(1-{\frac{(1-q)}{(N+2-Nq)}(x-\mu_q)^T\Sigma_q^{-1}(x-\mu_q)}
\right)_+^{\frac{1}{1-q}}
\nonumber
\end{align}
for all {$x\in\mathbb{R}^N$}, where 
{$\mu_q$} and {$\Sigma_q$} are known as the $q$-mean and $q$-covariance matrix, respectively.
These are generalizations of the usual mean and covariance~\cite{Prato_1999_jour_PhyRevE} and
correspond to the first and second moments with respect to the
so-called `deformed' expectation or normalized $q$-expectation that in turn is defined by
\begin{equation}
\label{q-expect-defn}
\langle{f}\rangle_q = 
\frac{\displaystyle\int_{\mathbb{R}^N}f(x) p(x)^q\:\mathrm{d}x}{\displaystyle\int_{\mathbb{R}^N}p(x)^q\:\mathrm{d}x},
\end{equation}
which is the expectation with respect to an escort distribution
{$p_q(x) = \frac{p(x)^q}{\int_{\mathbb{R}^N}p(x)^q\:\mathrm{d}x}$},
that is compatible with the foundations of nonextensive 
information theory~\cite{Tsallis_1998_jour_PhysicaA}.
The condition {$y_+=\max(y,0)$} in~\eqref{Gq:formula_q}, 
called the Tsallis cut-off condition~\cite{Tsallis_1995_jour_Chaos}, 
ensures that the above expression is well-defined, 
and {$K_{q,N}$} is the normalizing constant given by
\begin{equation}
K_{q,N} = \left\{\begin{array}{l}
		\left(\frac{N+2-Nq}{1-q}\right)^{\frac{N}{2}}
		\frac{\pi^{N/2} \Gamma\left(\frac{2-q}{1-q}\right)}{\Gamma\left(\frac{2-q}{1-q}+\frac{N}{2}\right)} 
		\hfill\text{~~~~for }  q<1,
		\\\\
		\left(\frac{N+2-Nq}{q-1}\right)^{\frac{N}{2}}
		\frac{\pi^{N/2} \Gamma\left(\frac{1}{q-1}-\frac{N}{2}\right)}{\Gamma\left(\frac{1}{q-1}\right)} 
		\hspace{15mm}
		\\ \hfill \text{for }  1<q<\left(1+\frac{2}{N}\right).
                \end{array}\right.
\label{eq_normalizing_const}
\end{equation}
The distribution~\eqref{Gq:formula_q} is only defined for $q<1+\frac{2}{N}$.
It retrieves the Gaussian distribution as $q\to1$,
and for $q>1$, it has a one-one correspondence with the
Student-$t$ distribution, with the special case of $q=1+\frac{2}{N+1}$ being the Cauchy distribution.
The uniform distribution on an infinitesimally small hypercube around the origin can
be obtained in the limit as $q\to-\infty$. 
For $q=0$, we have the smoothing kernel corresponding to random search algorithms~\cite{Rubinstein_1981_book_Wiley}.
Further, in the one-dimensional case, 
$q$-Gaussian with $q=-1$ gives the semicircle distribution, and
$q=\frac{\alpha-2}{\alpha-1}$ 
corresponds to the symmetric Beta$(\alpha,\alpha)$ distribution used in~\cite{Kreimer_1992_jour_AnnOR}.
In fact, the support of the 
$q$-Gaussian distribution can be expressed as
\begin{equation}
\Omega_q = \left\{\begin{array}{l}
		\left\{x: (x-\mu_q)^T\Sigma_q^{-1}(x-\mu_q) < \frac{N+2-Nq}{1-q}\right\}
		\\ \hfill\text{for }  q<1,
		\\\\
		\mathbb{R}^N 
		 \hfill\text{for } 1<q<\left(1+\frac{2}{N}\right).
                \end{array}\right.
\label{eq_support}
\end{equation}

Ghoshdastidar et al.~\cite{Ghoshdastidar_2013_arxiv} showed that the 
$q$-Gaussian family of distributions satisfy the Rubinstein 
conditions~\cite[pg 263]{Rubinstein_1981_book_Wiley} for smoothing kernels. The significance of the
$q$-Gaussian smoothing kernel is enhanced by the fact that it encompasses
the existing examples of smoothing kernels,
and thus significantly enhances the class of perturbation distributions for smoothed 
functional algorithms.
For the remainder of the paper, we will use an $N$-dimensional $q$-Gaussian distribution 
with zero $q$-mean and $q$-covariance matrix $\beta^2I_{N\times N}$. For convenience,
we refer to this distribution as $G_{q,\beta}(\cdot)$, with the 
case of $\beta=1$ being denoted by $G_q$.
We also use $\Omega_q$ to denote only the support set of $G_q$,
while we use $\theta+\beta\Omega_q$ for the support set of the distribution with $q$-mean
$\theta$ and $q$-covariance matrix $\beta^2I_{N\times N}$.
However, for $q>1$, the above set is always equal to $\mathbb{R}^N$.

A projected gradient based technique is commonly employed to optimize the smoothed
functional, and it has been observed 
in~\cite{Bhatnagar_2007_jour_TOMACS,Ghoshdastidar_2013_arxiv} 
that the two-sided gradient SF estimate provides significantly improved performance 
over the corresponding one-sided counterpart. 
In the case of $q$-Gaussian smoothing,
the gradient estimate is given by~\cite{Ghoshdastidar_2013_arxiv}
\begin{align}
\nabla_{\theta} S_{q,\beta}&[J(\theta)] 
= \mathsf{E}_{G_{q}(\eta)}
\left[\left.\frac{\eta \big(J(\theta+\beta\eta)-J(\theta-\beta\eta)\big)}
{\beta(N+2-Nq)\rho{(\eta)}} \right|\theta\right],
\label{grad_SF2_formula}
\end{align}
where 
the term
\begin{equation}
\rho(\eta) = \left(1-\frac{(1-q)}{(N+2-Nq)}\Vert{\eta}(n)\Vert^2\right)
\label{rho_defn}
\end{equation}
appears due to the differentiation of $G_q$.
It is shown that for $\beta$ small enough, the smoothed gradient is close to the
gradient of the objective function, assuming that it exists. 
Then,
a simple technique to estimate the above gradient is to consider a sample average
over some $L$ samples as
\begin{align}
\nabla_{\theta}&J(\theta) \approx \frac{1}{\beta L(N+2-Nq)} \times
\nonumber \\
&\sum_{n=0}^{L-1}
\frac{{\eta}(n)\big(J(\theta+\beta\eta(n))-J(\theta-\beta\eta(n))\big)}{\left(1-\frac{(1-q)}{(N+2-Nq)}\Vert{\eta}(n)\Vert^2\right)}\;.
\label{estimate1}
\end{align}

\subsection{Two-simulation $q$-Gaussian SF Hessian estimate}
In this section, we extend the above idea to the case of Hessians,
which is required in Newton based search algorithms.
Before presenting the estimate, we make a technical assumption
that ensures the existence of the gradient and Hessian of the objective function.
\begin{assumption}
\label{differentiable}
The function {$J(.)$} is twice continuously differentiable for all {$\theta\in C$}.
\end{assumption}

The above assumption is required for the theoretical analysis, but is
not necessary from a practical perspective,
since we hold $\beta>0$ fixed in the algorithm.
As with the case of gradient, 
the existence of {$\nabla_{\theta}^2 J(\theta)$} is assumed (Assumption~\ref{differentiable}), and
we estimate the same using SF approach. 
We define the two-sided smoothed Hessian with $q$-Gaussian smoothing
by following~\cite{Rubinstein_1981_book_Wiley}. 
For this, we can write the two-sided SF~\eqref{eq_2SF_general} as
\begin{align*}
S_{q,\beta} [J(\theta)] 
&= \frac{1}{2} \int\limits_{\beta\Omega_q}  G_{q,\beta}(\eta) J(\theta+\eta) \mathrm{d}\eta 
 \\&\qquad\qquad\nonumber
 +\frac{1}{2} \int\limits_{\beta\Omega_q}  G_{q,\beta}(\eta) J(\theta-\eta) \mathrm{d}\eta \;.
\end{align*}
Denoting the integrals by $\mathcal{S}_1(\theta)$ and $\mathcal{S}_2(\theta)$, respectively, and
substituting $\eta' = \theta+\eta$ in $\mathcal{S}_1(\theta)$, we have
\begin{align*}
 \mathcal{S}_1(\theta) 
 &= \frac{1}{\beta^N} \int\limits_{\theta+\beta\Omega_q}
 G_{q} \left(\frac{\eta'-\theta}{\beta}\right)J(\eta')\mathrm{d}\eta'
\end{align*}
where we use the fact that $G_{q,\beta}(\eta) = \frac{1}{\beta^N} G_q(\frac{\eta}{\beta})$,
which is true for all smoothing kernels.
The Hessian of $\mathcal{S}_1(\theta)$ is 
\begin{align*}
 \nabla_{\theta}^{2} \mathcal{S}_1(\theta)
 &= \frac{1}{\beta^N} \int\limits_{\theta+\beta\Omega_q}
 \nabla_{\theta}^{2} G_{q} \left(\frac{\eta'-\theta}{\beta}\right)J(\eta')\mathrm{d}\eta'\;.
\end{align*}
One can note that for $q<1$, the region over which integration is performed is a 
function of $\theta$, and hence by Leibnitz integral rule, there should be an additional
integral term, where the integration is over the surface of
the set $\theta+\beta\Omega_q$. However, since this integrand involves 
$G_{q} (\frac{\eta'-\theta}{\beta})$ that is zero over the surface, we can 
ignore the term completely.
Now, we substitute $\eta'' = \frac{\eta'-\theta}{\beta}$ above, and as a result, we have
$\mathrm{d}\eta'' = \frac{1}{\beta^N}\mathrm{d}\eta'$ and for all components
$i,j=1,\ldots,N$, $\frac{\partial{d}\eta''^{(i)}}{\partial{d}\theta^{(j)}} = -\frac{1}{\beta}$
whenever $i=j$, and 0 otherwise. Under this change of variables,
we can write
\begin{align}
 \nabla_{\theta}^{2} \mathcal{S}_1(\theta)
 &= \frac{1}{\beta^2} \int\limits_{\Omega_q}
 \nabla_{\eta''}^{2} G_{q} (\eta'') J(\theta+\beta\eta'')\mathrm{d}\eta''\;.
\end{align}
Similarly $\nabla_{\theta}^{2} \mathcal{S}_2(\theta)$ can also be derived,
and the two-sided smoothed Hessian with $q$-Gaussian smoothing is
\begin{align}
\nabla_{\theta}^2 &S_{q,\beta}[J(\theta)] \nonumber\\
  &= \frac{1}{2\beta^2} \int_{\Omega_q}  \nabla_{\eta}^{2} G_{q} (\eta) \big(J(\theta+\beta\eta) 
  + J(\theta-\beta\eta) \big) \mathrm{d}\eta \;.
\label{hess_SF_defn}
\end{align}

We now compute the Hessian matrix corresponding to the standard $q$-Gaussian distribution, $G_q$.
When $\eta\in\Omega_q$,
the partial derivative of $G_q(\eta)$ with respect to $\eta^{(i)}$ for all $i=1,\ldots,N$ is given by
\begin{align*}
 \frac{\partial G_q(\eta)}{\partial \eta^{(i)}} 
 = -\frac{2\eta^{(i)} \rho(n)^{\frac{q}{1-q}}}{K_{q,N}(N+2-Nq)}\;,
\end{align*}
where $\rho(\cdot)$ is as defined in~\eqref{rho_defn}.
From above, we can compute the second derivatives,
which can be expressed in terms of $G_q(\eta)$ as
\begin{align*}
\frac{\partial^2 G_q(\eta)}{\partial \eta^{(i)} \partial \eta^{(j)}}
= \frac{4q \eta^{(i)} \eta^{(j)}}{(N+2-Nq)^2} \frac{G_q(\eta)}{\rho(\eta)^2}
\end{align*}
for {$i\neq j$}. For {$i=j$}, we have
\begin{align*}
\frac{\partial^2 G_q(\eta)}{\partial {\eta^{(i)}}^2}
&= \frac{4q {\eta^{(i)}}^2}{(N+2-Nq)^2} \frac{G_q(\eta)}{\rho(\eta)^2}
\\&\qquad\qquad
- \frac{2}{(N+2-Nq)} \frac{G_q(\eta)}{\rho(\eta)} \;.
\end{align*}

Thus, the Hessian turns out to be of the form
\begin{align}
\nabla_{\eta}^2 G_{q}(\eta)
= \frac{2}{(N+2-Nq)} H(\eta) G_{q}(\eta),
\end{align}
where
\begin{align}
H(\eta) = \left\{ \begin{array}{rcl} 
          \left(\displaystyle\frac{2q}{(N+2-Nq)} \displaystyle\frac{\eta^{(i)}\eta^{(j)}}
          {\rho{(\eta)}^2}\right) &\text{for} &i\neq j
          \\ \\
          \left(\displaystyle\frac{2q}{(N+2-Nq)} \displaystyle\frac{\left(\eta^{(i)}\right)^2}{\rho{(\eta)}^2} 
           - \displaystyle\frac{1}{\rho{(\eta)}}\right) &\text{for} &i=j.
          \end{array} \right. 
\label{H_defn}
\end{align}
is a generalization of a similar function given in~\cite{Bhatnagar_2007_jour_TOMACS}, 
that can be obtained as {$q\to1$}. 
%
Substituting $\nabla_{\eta}^2 G_{q}(\eta)$ in~\eqref{hess_SF_defn}, we have
\begin{align}
&\nabla_{\theta}^2 S_{q,\beta}[J(\theta)] 
\nonumber \\
&=
\mathsf{E}_{G_{q}(\eta)}\left[ \left. \frac{H(\eta) \big(J(\theta+\beta\eta)+J(\theta-\beta\eta)\big)}{\beta^2(N+2-Nq)} \right| \theta \right].
\label{hess_SF2_formula}
\end{align}
Subsequently, we show that the Hessians of $S_{q,\beta}[J(\theta)]$ and $J(\theta)$
are close enough, and hence,
 we obtain the Hessian  estimate of $J(\theta)$, for large $L$ and small {$\beta$}, as 
\begin{align}
\nabla_{\theta}^2 &J(\theta) \approx\frac{1}{\beta^2 L (N+2-Nq)} \times
\nonumber \\
&\sum_{n=0}^{L-1}
H\big(\eta(n)\big) \big(J(\theta+\beta\eta(n))+J(\theta-\beta\eta(n))\big) .
\label{SF_estimate1}
\end{align}
In the next section, we present a Newton based search technique
using the above gradient and Hessian estimates,~\eqref{estimate1} 
and~\eqref{SF_estimate1}, respectively.

\section{Optimization of long-run average cost of a parametrized  Markov process}
\label{sec_problem}

Our objective here is to optimize $J:C\mapsto\mathbb{R}$,
when only (noisy) observations of $J$ are known.
If $J$ is an analytic function, as considered 
in~\cite{Spall_2000_jour_AutoCtrl,Chin_1997_jour_TransSysManCyber}, then the estimates
in~\eqref{estimate1} and~\eqref{SF_estimate1} can be directly used.
We consider a slightly complicated scenario,
often encountered in problems of stochastic control,
where the objective is a stochastic function with no analytic expression.
Such a setting is discussed below.

\subsection{Problem Framework}

Let {$\{Y_n\}_{n\geqslant0} \subset\mathbb{R}^d$} be a parameterized Markov process
with transition kernel {${P}_{\theta}(x,\:\mathrm{d}y)$} that depends
on a tunable parameter {$\theta\in C$}, where {$C\subset\mathbb{R}^N$} is compact and convex.
We assume the following.
%
\begin{assumption}
\label{ergodic}
For a fixed operative parameter $\theta\in C$,
the Markov process {$\{Y_n\}$} is ergodic and has a unique invariant measure $\nu_{\theta}$. 
\end{assumption}
%
Though we restrict ourselves to an ergodic Markov process in this paper, the subsequent 
discussions can be directly extended to hidden Markov models following the lines 
of~\cite{Bhatnagar_2001_jour_IIETrans}. Thus the work in this paper is also applicable,
with suitable modifications, to a broader class of problems.
We also consider a Lipschitz continuous cost function 
{$h:\mathbb{R}^d\mapsto\mathbb{R}^+\bigcup \{0\}$} associated with the process.
Our objective is to minimize the long-run average cost
\begin{equation}
\label{Jdefn}
J(\theta) = \lim_{M\to\infty}\frac{1}{M}\sum_{m=0}^{M-1}h(Y_m) = \int\limits_{\mathbb{R}^d}h(x)\nu_{\theta}(\:\mathrm{d}x),
\end{equation}
by choosing an appropriate {$\theta\in C$}. The existence of the above limit is assured by Assumption~\ref{ergodic}
and the fact that $h$ is continuous, hence measurable. 
In addition, we assume that the average cost {$J({\theta})$} satisfies 
Assumption~\ref{differentiable}.
However, in this setting, verification of Assumption~\ref{differentiable} depends on the underlying process
and is non-trivial in most cases. One can observe that under certain conditions
(for instance when cost function $h(\cdot)$ is bounded),
Assumption~\ref{differentiable} can be translated to impose the condition of continuous
differentiability of the stationary measure $\nu_\theta$ for all $\theta\in C$.
This, in turn, would depend on a similar condition on the transition kernel $P_\theta(x,\mathrm{d}y)$.
Discussions on such conditions for finite state Markov processes can be found 
in~\cite{Schweitzer_1968_jour_AppProb,Kushner_1981_jour_SIAMCtrlOptim}, and similar 
results for general state systems were presented in~\cite{Vazquez_1992_jour_AppProb}.
However, in the general case, such conditions are difficult to verify.
In addition to above, we also assume the existence of a stochastic Lyapunov function. 
This requires the notion of a non-anticipative sequence, defined below.

\begin{definition}[Non-anticipative sequence]
Any random sequence of parameter vectors, {$(\theta(n))_{n\geqslant0} \subset C$}, controlling a process
{$\{Y_n\}\subset\mathbb{R}^d$}, is said to be non-anticipative if the conditional probability
{$P(Y_{n+1}\in B |\mathcal{F}_n) = {P}_{\theta} (Y_n, B)$} 
almost surely for {$n\geqslant0$} and all Borel sets {$B\subset\mathbb{R}^d$},
where {$\mathcal{F}_n = \sigma(\theta(m),Y_m,m\leqslant n)$}, {$n\geqslant0$} are the associated {$\sigma$}-fields. 
\end{definition}
One can verify that under a non-anticipative parameter sequence {$(\theta(n))$}, the joint process
{$(Y_n,\theta(n))_{n\geqslant0}$} is Markov.
We assume the existence of a stochastic Lyapunov function (below), which
ensures that the process under a tunable parameter remains stable.
Assumption~\ref{lyapunov} will not be required, for instance, if, in addition, the single-stage cost function $h$ is bounded. 
It can be seen that the sequence of parameters obtained using our algorithm
forms a non-anticipative sequence.

\begin{assumption}
\label{lyapunov}
Let {$(\theta(n))$} be a non-anticipative sequence of random parameters controlling the process {$\{Y_n\}$},
and {$\mathcal{F}_n = \sigma(\theta(m),$} {$Y_m,m\leqslant n)$}, {$n\geqslant0$} be a sequence of associated {$\sigma$}-fields. 
There exists {$\epsilon_0>0$}, a compact set {$\mathcal{K}\subset\mathbb{R}^d$}, and a continuous function 
{$V:\mathbb{R}^d\mapsto\mathbb{R}^+\bigcup\{0\}$}, with
{$\lim_{\Vert{x}\Vert\to\infty} V(x) = \infty$}, such that 
\begin{enumerate}
\item 
{$\sup\limits_n \mathsf{E}[V(Y_n)^2] < \infty$}, and
\item 
{$\mathsf{E}[V(Y_{n+1})|\mathcal{F}_n] \leqslant V(Y_n) - \epsilon_0$}, whenever {$Y_n\notin\mathcal{K}$}, {$n\geqslant0$}.
\end{enumerate} 
\end{assumption}

As a consequence of Assumption~\ref{ergodic},
we can estimate the gradient and Hessian,~\eqref{estimate1} and~\eqref{SF_estimate1}
respectively, as
\begin{align}
\nabla_{\theta}J(\theta) &\approx\frac{1}{\beta ML(N+2-Nq)} \times
\nonumber \\
&\sum_{n=0}^{M-1}\sum_{m=0}^{L-1}
\frac{{\eta}(n)\big(h(Y_m)-h(Y'_m)\big)}{\left(1-\frac{(1-q)}{(N+2-Nq)}\Vert{\eta}(n)\Vert^2\right)}
\label{estimate2}
\end{align}
and
\begin{align}
\nabla_{\theta}^2 J(\theta) &\approx\frac{1}{\beta^2 ML (N+2-Nq)} \times
\nonumber \\
&\sum_{n=0}^{M-1}\sum_{m=0}^{L-1} 
H(\eta(n)) \big(h(Y_m)+h(Y'_m)\big),
\label{SF_estimate2}
\end{align}
for large $M$, $L$ and small {$\beta$}, where {$\{Y_m\}$} and {$\{Y'_m\}$} are governed by the parameters
{$(\theta+\beta{\eta}(n))$} and {$(\theta-\beta{\eta}(n))$}, respectively.

\subsection{Proposed Newton based technique}

Since $C$ is a compact and convex
subset of $\mathbb{R}^N$, projected gradient or Newton methods can be used,
where the update rule is of the form
\begin{equation}
\theta(n) = \mathcal{P}_C \big(\theta(n-1) - a(n)Z(n)\big)
\label{eq_graddes_update}
\end{equation}
for gradient based search, and
\begin{equation}
\theta(n) = \mathcal{P}_C \left(\theta(n-1) - a(n)W(n)^{-1}Z(n)\right)
\label{eq_newton_update}
\end{equation}
for Newton's method. Here, $\mathcal{P}_C$ is a projection operator onto the set $C$,
$Z(n)$ and $W(n)$ are estimates of the gradient vector and the Hessian matrix, respectively,
of the objective function at the $n^{th}$ iteration, 
and $(a(n))_{n\geqslant0}$ is a prescribed non-increasing step-size sequence.
The update in~\eqref{eq_graddes_update} corresponds to the gradient based $q$-Gaussian
SF algorithms~\cite{Ghoshdastidar_2013_arxiv}.

The estimators for the gradient and Hessian, given in~\eqref{estimate2} and~\eqref{SF_estimate2}, respectively,
are quite computationally intensive considering the fact that, at each iteration,
we require the sample size to be considerably large so that 
the steady state average for a given parameter update can be approximated closely.
A computationally efficient solution to this 
problem is to consider a multi-timescale stochastic approximation 
scheme~\cite{Bhatnagar_1998_jour_ProbEnggInfoSc}.
The idea is to update the estimates, $Z$ and $W$, on a different timescale,
faster than the update iteration as
\begin{align}
Z(n+1) &= (1-b(n))Z(n) + b(n)\hat{Z}(n),
\label{eq_2time_Z}
\\
W(n+1) &= \mathcal{P}_{pd} \big((1-b(n))W(n) + b(n)\hat{W}(n)\big), 
\label{eq_2time_W}
\end{align}
where $Z(n), W(n)$ are the updates till the $n^{th}$ iteration,
and $\hat{Z}(n), \hat{W}(n)$ are the instantaneous estimates 
of the gradient and Hessian, 
using~\eqref{estimate2} and~\eqref{SF_estimate2} where one may let $M=1$.
Also $\mathcal{P}_{pd}$ is an operator that projects any $N \times N$ matrix to
the space of positive definite and symmetric matrices.
We make the following 
assumption on the two step-size sequences
{$(a(n))_{n\geqslant0}$} and {$(b(n))_{n\geqslant0}$}.

\begin{assumption}
\label{stepsize}
{$(a(n))_{n\geqslant0}$}, {$(b(n))_{n\geqslant0}$} are sequences 
of positive scalars such  that
\begin{enumerate}
\item
$b(n) \leqslant 1$ for all $n$,
\item
{$a(n) = o(b(n))$}, \textit{i.e.}, {$\frac{a(n)}{b(n)}\to0$} as {$n\to\infty$},
\item
{$\sum_{n=0}^{\infty}a(n) = \sum_{n=0}^{\infty}b(n) = \infty$},
\item
{$\sum_{n=0}^{\infty}a(n)^2 <\infty$}, and {$\sum_{n=0}^{\infty}b(n)^2 <\infty$}.
\end{enumerate}
\end{assumption}

Considering the step-sizes as above, if we update both gradient and Hessian
estimate using the larger step-sizes $(b(n))$,
then these estimates are updated on a faster timescale.
On the other hand, the parameter $\theta(n)$, when updated using step-sizes $(a(n))$
appears to change slowly and is seen to be nearly constant from the timescale of $b(n)$.
Thus, even if we choose $M$ in~\eqref{estimate2} and~\eqref{SF_estimate2}
to be small (say $M=1$), the asymptotic stationarity of the process
is not affected because the updates of the estimators always occur with the parameter 
in a quasi-static state.
The convergence analysis in Section~\ref{sec_convergence} shows that
when estimated as per~\eqref{eq_2time_Z}--\eqref{eq_2time_W}, the
sample size $L$ can also be set as $L=1$. However, it has been
observed~\cite{Bhatnagar_2003_jour_Simulation} that updates along a subsample, 
typically $L$ between 50 and 100 gives desirable performance.

One of the issues with Newton-based algorithms is that the Hessian has to be positive definite
for the algorithm to progress in the descent direction. This 
may not hold always.
Hence, the estimate obtained at each update step has to be projected onto
the space of positive definite and symmetric matrices. 
This is taken care of by the map
{$\mathcal{P}_{pd}:\mathbb{R}^{N\times N}\mapsto\{$}symmetric matrices with eigenvalues{$\geqslant\varepsilon\}$}
that projects any {$N\times N$} matrix onto the 
set of symmetric positive definite matrices, with a minimum eigenvalue 
of at least {$\varepsilon$} for some {$\varepsilon>0$}. 
We assume the projection operator $\mathcal{P}_{pd}$ satisfies the following:
\begin{assumption}
\label{projection}
If {$(A_n)_{n\in\mathbb{N}}$}, {$(B_n)_{n\in\mathbb{N}} \subset \mathbb{R}^{N\times N}$} are sequences of matrices
satisfying {$\lim_{n\to\infty} \Vert{A_n-B_n}\Vert = 0$}, then 
{$\lim_{n\to\infty} \Vert{\mathcal{P}_{pd}(A_n)-\mathcal{P}_{pd}(B_n)}\Vert = 0$}. 
\end{assumption}

We present our Newton based algorithm below, which we denote N$q$-SF2 to signify that
it employs a Newton based approach with two-simulation $q$-Gaussian SF.
The update runs for some specified $M$ iterations 
(not to be confused with the $M$ in the estimators above), and
$Z$ and $W$ denote the estimators for the gradient and the Hessian, respectively.
It requires sampling from a standard $q$-Gaussian
distribution for some $q\in(0,1+\frac{2}{N})$. The procedure for generating $q$-Gaussian samples is
given in~\cite{Ghoshdastidar_2013_arxiv}. The reason behind the restriction on
the values of $q$ is discussed in the next section. 

\begin{alg}[N$q$-SF2 algorithm]
\label{NqSF2_problem2}
Assuming that constants $\beta>0$, $\varepsilon>0$, 
$q\in(0,1+\frac{2}{N})$, $L$ and $M$,
and the step-sizes {$(a(n))_{n\geqslant0}$}, {$(b(n))_{n\geqslant0}$}
are specified, the algorithm proceeds as below.
\begin{enumerate}
\item
Initialize some $\theta(0) \in C$.
\item
Set $Z=0 \in \mathbb{R}^N$ and {$W(0) = 0 \in\mathbb{R}^{N\times N}$}.
\item
For $n=0$ to $M-1$
	\begin{enumerate}
	\item
	Generate a random vector {$\eta(n) \in \mathbb{R}^N$}
	from a standard $N$-dimensional $q$-Gaussian distribution\;
	\item
	For $m=0$ to $L-1$
		\begin{enumerate}
		\item
		Generate two independent simulations {$Y_{nL+m}$} and
		{$Y'_{nL+m}$} governed by 
		{$\mathcal{P}_C(\theta(n)+\beta\eta(n))$} and 
		{$\mathcal{P}_C(\theta(n)-\beta\eta(n))$}, respectively.
		\item
		Update gradient estimate as
		\begin{align*}
		& \hspace{-10mm}
		Z(nL+m+1) = (1-b(n))Z(nL+m) +\\
		& \hspace{-9mm}
		b(n) \left[\frac{\eta(n)(h(Y_{nL+m})-h(Y'_{nL+m}))}
		{\beta(N+2-Nq)(1-\frac{(1-q)}{(N+2-Nq)}\Vert{\eta}(n)\Vert^2)}\right]
		\end{align*}
		\item
		Compute $H(\eta(n))$ using~\eqref{H_defn}, and 
		update Hessian estimate as
		\begin{align*}
		& \hspace{-10mm}
		W(nL+m+1) = (1-b(n))W(nL+m) +\\
		& \hspace{-9mm}
		b(n) \left[\frac{H(\eta(n))(h(Y_{nL+m})+h(Y'_{nL+m}))}
		{\beta^2(N+2-Nq)}\right]
		\end{align*}
		\end{enumerate}
	\item
	Project Hessian matrix, \textit{i.e.}, 
	\begin{displaymath}
	W((n+1)L) := \mathcal{P}_{pd}\big(W((n+1)L)\big).
	\end{displaymath}
	\item
	Update $\theta(n+1) = $
	\begin{displaymath}
	\hspace{-7mm}
	\mathcal{P}_C \left( \theta(n) - a(n)W((n+1)L)^{-1}Z((n+1)L) \right).
	\end{displaymath}
	\end{enumerate}
\item
Output $\theta(M)$ as the final parameters.
\end{enumerate}
\end{alg}

For implementation purposes, a modified version of the above algorithm is often 
considered~\cite{Spall_2000_jour_AutoCtrl,Bhatnagar_2007_jour_TOMACS}.
This is known as the Jacobi variant, where the projection map $\mathcal{P}_{pd}$
is such that it sets all the off-diagonal terms in $W(n)$ to zero, and the diagonal
terms are projected onto the interval $[\varepsilon,\infty)$. 
This ensures that the projected matrix has a minimum eigenvalue of at least $\varepsilon$,
and this also simplifies the inverse computation.

\section{Convergence of the proposed algorithm}
\label{sec_convergence}

We present below our main convergence results whose proofs can be 
found in Appendix~\ref{app_convergence} at the end of the paper.
Let us consider the updates along the faster timescale, \textit{i.e.}, Step~(3b) of the N$q$-SF2 algorithm.

We define {$\tilde{\theta}(p) = \theta(n)$}, {$\tilde{\eta}(p) = \eta(n)$} and
{$\tilde{b}(p) = b(n)$} for {$nL\leqslant p<(n+1)L$}, {$n\geqslant0$}. From Assumption~\ref{stepsize}, we have
{$a(p) = o\big(\tilde{b}(p)\big)$}, {$\sum_p \tilde{b}(p) = \infty$} and {$\sum_p \tilde{b}(p)^2 < \infty$}. 
Since, {$\{Y_p\}_{p\in\mathbb{N}}$} and {$\{Y'_p\}_{p\in\mathbb{N}}$} are independent Markov processes,
we can consider {$\{(Y_p,Y'_p)\}_{p\in\mathbb{N}}$} as a joint Markov process 
parameterized by 
{$\big(\mathcal{P}_C{(\tilde{\theta}(p)+\beta\tilde{\eta}(p))},\mathcal{P}_C{(\tilde{\theta}(p)-\beta\tilde{\eta}(p))}\big)$}.
We can rewrite Step~(3b.ii) using the following iteration for all {$p\geqslant0$}:
\begin{equation}
Z(p+1) = Z(p) + \tilde{b}(p)\big( g_1(Y_p,Y'_p,\tilde{\eta}(p)) - Z(p) \big),
\label{faststep_Z}
\end{equation}
where  
\begin{equation}
 g_1(Y_p,Y'_p,\tilde{\eta}(p)) = 
 \left(\displaystyle\frac{\tilde{\eta}(p) (h(Y_p)-h(Y'_p))}{\beta (N+2-Nq) \rho{(\tilde{\eta}(p))}}\right)
\end{equation}
for {$nL\leqslant p<(n+1)L$}, with {$\rho{(.)}$} defined as in~\eqref{rho_defn}.
Similarly, the update of the Hessian matrix in Step~(3b.iii) can be expressed as
\begin{equation}
W(p+1) = W(p) + \tilde{b}(p)\big( g_2(Y_p,Y'_p,\tilde{\eta}(p)) - W(p) \big),
\label{faststep_W}
\end{equation}
where, for {$nL\leqslant p<(n+1)L$}, 
\begin{equation}
 g_2(Y_p,Y'_p,\tilde{\eta}(p)) = \left(\displaystyle\frac{{H}(\tilde{\eta}(p)) (h(Y_p)+h(Y'_p))}{\beta^2(N+2-Nq)}\right).
\end{equation}
Let {$\mathcal{G}_p = \sigma\big(\tilde{\theta}(k), \tilde{\eta}(k), Y_k, Y'_k, k\leqslant p\big), p\geqslant0$} 
denote a sequence of $\sigma$-fields
generated by the mentioned quantities. We can observe that {$(\mathcal{G}_p)_{p\geqslant0}$} is a filtration,
where {$g_1(Y_p,Y'_p,\tilde{\eta}(p))$} and {$g_2(Y_p,Y'_p,\tilde{\eta}(p))$} 
are {$\mathcal{G}_p$}-measurable for each {$p\geqslant0$}. 
We can rewrite~\eqref{faststep_Z} and~\eqref{faststep_W} as
\begin{align}
Z&(p+1) = Z(p) + 
\nonumber \\&\tilde{b}(p)\big[ E[g_1(Y_p,Y'_p,\tilde{\eta}(p))|\mathcal{G}_{p-1}] - Z(p) + A_p\big],
\label{eq_Z_step}
\\
W&(p+1) = W(p) + 
\nonumber \\&\tilde{b}(p)\big[ E[g_2(Y_p,Y'_p,\tilde{\eta}(p))|\mathcal{G}_{p-1}] - W(p) + B_p\big],
\label{eq_W_step}
\end{align}
where {$A_p = g_1(Y_p,Y'_p,\tilde{\eta}(p)) - E[g_1(Y_p,Y'_p,\tilde{\eta}(p))|\mathcal{G}_{p-1}]$} and
{$B_p = g_2(Y_p,Y'_p,\tilde{\eta}(p)) - E[g_2(Y_p,Y'_p,\tilde{\eta}(p))|\mathcal{G}_{p-1}]$} 
are both {$\mathcal{G}_p$}-measurable. 

The following result presents a useful property of {$(A_p)_{p\in\mathbb{N}}$} 
and {$(B_p)_{p\in\mathbb{N}}$}. 

\begin{lem}
\label{B_n_convergence}
For all values of {$q\in\big(-\infty,1\big)\cup\big(1,1+\frac{2}{N}\big)$},
{$(A_p,\mathcal{G}_p)_{p\in\mathbb{N}}$} and {$(B_p,\mathcal{G}_p)_{p\in\mathbb{N}}$} 
are martingale difference sequences with bounded variance. 
\end{lem}

The iterations~\eqref{faststep_Z} and~\eqref{faststep_W} are 
not coupled, \textit{i.e.}, iterates $Z(p)$ do not depend on $W(p)$ 
and vice-versa. Thus,
they can be dealt with separately. 
We can write the parameter update along the slower timescale as
{$\theta(n+1) = \mathcal{P}_C\big(\theta(n) - \tilde{b}(n)\zeta(n)\big)$},
where we use 
\begin{displaymath}
 \zeta(n)=\frac{a(n)}{\tilde{b}(n)}W((n+1)L)^{-1}Z((n+1)L)=o(1),
\end{displaymath}
since {$a(n) = o(\tilde{b}(n))$}. 
Thus, the parameter update recursion is quasi-static
when viewed from the timescale of {$(\tilde{b}(n))$}, and hence, one may let
{$\tilde{\theta}(p) \equiv \theta$} and {$\tilde{\eta}(p) \equiv \eta$} for all {$p\in\mathbb{N}$}, when analyzing~\eqref{eq_Z_step}
and \eqref{eq_W_step}.
The system of ODEs associated with these updates is the following:
\begin{align}
\dot{\theta}(t) &= 0,
\label{fastode_2}
\\ \dot{Z}(t) &= \frac{\eta \big(J(\theta+\beta\eta)-J(\theta-\beta\eta)\big)}{\beta(N+2-Nq)\rho{(\eta)}}  - Z(t)\;,
\label{fastode_Z}
\\ \text{and~~}
 \dot{W}(t) &= \frac{{H}(\eta) \big(J(\theta+\beta\eta)+J(\theta-\beta\eta)\big)}{\beta^2(N+2-Nq)}  - W(t)\;.
\label{fastode_W}
\end{align}

At this stage, we recall a series of 
results by Borkar~\cite{Borkar_2008_book_Cambridge}.
\begin{thm}
\label{borkar_cor_8}
\textbf{\cite[Thm 7--Cor 8, pp. 74 and Thm~9, pp. 75]{Borkar_2008_book_Cambridge}}
Consider the iteration, 
\begin{displaymath}
 x_{p+1} = x_p + \gamma(p)\big[f(x_p,Y_p)+M_p\big].
\end{displaymath}
Let the following conditions hold:
\begin{enumerate}
\item 
{$\{Y_p:p\in\mathbb{N}\}$} is a Markov process satisfying Assumptions~\ref{ergodic} and~\ref{lyapunov},
\item
for each {$x\in\mathbb{R}^N$} and {$x_p\equiv x$} for all $p\in\mathbb{N}$,
{$Y_p$} has a unique invariant probability measure {$\nu_{x}$},
\item
{$(\gamma(p))_{p\geqslant0}$} are step-sizes satisfying {$\sum\limits_{p=0}^{\infty} \gamma(p)=\infty$} 
and {$\sum\limits_{p=0}^{\infty} \gamma^2(p)<\infty$},
\item
{$f(.,.)$} is Lipschitz continuous in its first argument uniformly w.r.t the second, 
\item
{$M_p$} is a martingale difference noise term with bounded variance,
\item
if {$\tilde{f}\big(x,\nu_x\big) = \mathsf{E}_{\nu_x} \big[f(x,Y)\big]$}, then the limit
\begin{displaymath}
 \hat{f}\big(x(t)\big) = \displaystyle\lim\limits_{a\uparrow\infty} \frac{\tilde{f}\big(ax(t),\nu_{ax(t)}\big)}{a}
\end{displaymath}
exists uniformly on compacts, and
\item
the ODE {$\dot{x}(t) = \hat{f}\big(x(t)\big)$} is well-posed
and has the origin as the unique globally asymptotically stable equilibrium.
\end{enumerate}
Then the update {$x_p$} satisfies {$\sup_p \Vert{x_p}\Vert < \infty$}, almost surely,
and converges to the stable fixed points of the ODE
\begin{displaymath}
\dot{x}(t) = \tilde{f}\big(x(t),\nu_{x(t)}\big).
\end{displaymath}
\end{thm}

As a consequence of Lemma~\ref{B_n_convergence} and the above result,
we have the following lemma proving the convergence of the gradient and Hessian updates.
\begin{lem}
\label{W_bounded}
The sequences $(Z(p))$ and {$(W(p))$} are uniformly bounded with probability 1. Further,
\begin{align*}
 &\left\Vert  Z(p) - 
 \frac{\tilde{\eta}(p)\big(J(\tilde{\theta}(p)+\beta\tilde{\eta}(p))-J(\tilde{\theta}(p)-\beta\tilde{\eta}(p))\big)}{
 \beta(N+2-Nq)\rho{(\tilde{\eta}(p))}}
 \right\Vert,
\\
 &\left\Vert W(p) -
 \frac{{H}(\tilde{\eta}(p))\big(J(\tilde{\theta}(p)+\beta\tilde{\eta}(p))+J(\tilde{\theta}(p)-\beta\tilde{\eta}(p))\big)}{
 \beta^2 (N+2-Nq)} 
 \right\Vert
\end{align*}
$\to 0$ almost surely as {$p\to\infty$}. 
\end{lem}

Thus, both {$Z$} and {$W$} recursions eventually track the gradient and Hessian of {$S_{q,\beta}[J(\theta)]$}.
So, after incorporating the projection considered in Step~(3c), we can write the parameter update, Step~(3d)
of the N$q$-SF2 algorithm, as 
\begin{align}
\theta (n+&1)= \mathcal{P}_C \bigg(\theta(n) + a(n) \bigg[ \Delta \big(\theta(n)\big) + \xi_n
\nonumber \\
& -\mathcal{P}_{pd}\left(\nabla_{\theta(n)}^2 J\big(\theta(n)\big)\right)^{-1}
\nabla_{\theta(n)} J\big(\theta(n)\big) \bigg]\bigg),
\label{hess_slowstep}
\end{align}
where we use~\eqref{grad_SF2_formula} and~\eqref{hess_SF2_formula} to write
\begin{align}
&\Delta \big(\theta(n)\big) =
\mathcal{P}_{pd}\left(\nabla_{\theta(n)}^2 J\big(\theta(n)\big)\right)^{-1}
\nabla_{\theta(n)} J\big(\theta(n)\big)
\nonumber \\&
- \mathcal{P}_{pd}\left(\nabla_{\theta(n)}^2 S_{q,\beta}\Big[J\big(\theta(n)\big)\Big]\right)^{-1}
\nabla_{\theta(n)} S_{q,\beta}\Big[J\big(\theta(n)\big)\Big],
\label{hess_error_term}
\end{align}
and the noise term 
\begin{align}
&\xi_n =
\mathsf{E}\bigg[\mathcal{P}_{pd}
\left(\frac{H(\eta(n))\bar{J}_n}{\beta^2 (N+2-Nq)}\right)^{-1} \times
\nonumber\\&\qquad\qquad\qquad\qquad\qquad
\frac{\eta(n)\bar{J}'_n}{\beta\rho{(\eta(n))}(N+2-Nq)}\bigg|\theta(n)\bigg]
\nonumber \\&
- \mathcal{P}_{pd}
\left(\frac{H(\eta(n))\bar{J}_n}{\beta^2 (N+2-Nq)}\right)^{-1}
\frac{\eta(n)\bar{J}'_n}{\beta\rho{(\eta(n))}(N+2-Nq)}\;,
\label{hess_noise_term}
\end{align}
where {$\bar{J}_n= J\big(\theta(n)+\beta\eta(n)\big)+J\big(\theta(n)-\beta\eta(n)\big)$}
and {$\bar{J}'_n= J\big(\theta(n)+\beta\eta(n)\big)-J\big(\theta(n)-\beta\eta(n)\big)$}.
It may be noted that the second term in~\eqref{hess_error_term} is the same as the first in~\eqref{hess_noise_term}. 
The next few results discuss some properties of  
the error term {$\Delta \big(\theta(n)\big)$} and the noise term {$\xi_n$},
that will be used to prove the convergence of N$q$-SF2 to a local optimum.

\begin{prop}
\label{hess_SF2_convergence}
For a given {$q\in\big(0,1\big)\bigcup\big(1,1+\frac{2}{N}\big)$}, for all {$\theta\in C$} and {$\beta>0$},
\begin{align*}
& \left\Vert \nabla_{\theta} S_{q,\beta}[J(\theta)]-\nabla_{\theta} J(\theta) \right\Vert = o(\beta)
\\ \text{and }&
\left\Vert \nabla_{\theta}^2 S_{q,\beta}[J(\theta)]-\nabla_{\theta}^2 J(\theta) \right\Vert = o(\beta). 
\end{align*}
Further, if Assumption~\ref{projection} holds, then {$\left\Vert\Delta (\theta)\right\Vert = o(\beta)$}. 
\end{prop}

One may note that the proof of Proposition~\ref{hess_SF2_convergence} imposes the condition of $q>0$. Hence, subsequent
analysis and simulations of N$q$-SF2  algorithm have been done only over the range of $q$ specified above.
The following result deals with the noise term $\xi_n$.
For this we consider the filtration {$(\mathcal{F}_n)_{n\geqslant0}$} defined as
{$\mathcal{F}_n=\sigma\big(\theta(0),\ldots,\theta(n),\eta(0),\ldots,\eta(n-1)\big)$}. 

\begin{lem}
\label{hess_xi_n_convergence}
Defining {$M_n = \sum_{i=0}^{n-1} a(k)\xi_k$}, {$(M_n,\mathcal{F}_n)_{n\geqslant0}$} 
is an almost surely convergent martingale sequence for all  {$q\in\big(0,1\big)\bigcup\big(1,1+\frac{2}{N}\big)$}.
\end{lem}
We state the following result due to Kushner and Clark~\cite{Kushner_1978_book_Springer} adapted to our scenario.
\begin{lem}
\label{kushner_thm_5.3.1}
\textbf{\cite[Theorem~5.3.1, pp 189--196]{Kushner_1978_book_Springer}}
Given the iteration, {$x_{n+1} = \mathcal{P}_C \big(x_n + \gamma_n(f(x_n) + \xi_n)\big)$}, where
\begin{enumerate}
\item
{$\mathcal{P}_C$} represents a projection operator onto a closed and bounded constraint set {$C$}, 
\item 
{$f(.)$} is a continuous function,
\item
{$(\gamma_n)_{n\geqslant0}$} is a positive sequence satisfying {$\gamma_n\downarrow0$}, {$\sum_{n=0}^{\infty} \gamma_n=\infty$}, and 
\item
{$\sum_{n=0}^m \gamma_n\xi_n$} converges a.s.
\end{enumerate}
Under the above conditions, the update {$(x_n)$} converges almost surely 
to the set of asymptotically stable fixed points of the ODE
\begin{equation}
\label{projectedode}
\dot{x}(t) = \tilde{\mathcal{P}}_C \big(f(x(t))\big),
\end{equation}
where {$\tilde{\mathcal{P}}_C \big(f(x)\big) = 
\lim\limits_{\epsilon\downarrow0}\left(\frac{\mathcal{P}_C \big(x+\epsilon f(x)\big)-x}{\epsilon}\right)$}.
\end{lem}
Proposition~\ref{hess_SF2_convergence} and Lemma~\ref{hess_xi_n_convergence} 
can be combined with Lemma~\ref{kushner_thm_5.3.1}
to derive the main theorem which affirms the convergence of the N$q$-SF2 algorithm. 
\begin{thm}
\label{thm_NqSF2}
Under Assumptions \ref{differentiable} -- \ref{projection}, given {$\epsilon>0$} and {$q\in\big(0,1\big)\bigcup\big(1,1+\frac{2}{N}\big)$},
there exists {$\beta_0 >0$} such that
for all {$\beta\in(0,\beta_0]$}, the sequence {$(\theta(n))$} obtained using N$q$-SF2 converges almost surely
as {$n\to\infty$} to the $\epsilon$-neighborhood of the set of stable attractors of the ODE
\begin{equation}
\dot{\theta}(t) = \tilde{\mathcal{P}}_C \left(
\mathcal{P}_{pd}\left(\nabla_{\theta(t)}^2 J\big(\theta(t)\big)\right)^{-1}
\nabla_{\theta(t)} J\big(\theta(t)\big)\right)
\label{New_ode}
\end{equation}
where the domain of attraction is 
\begin{equation}
\left\{\theta\in C \left| \nabla_{\theta}J(\theta)^T
\tilde{\mathcal{P}}_C \left(-
\mathcal{P}_{pd}\left(\nabla_{\theta}^2 J(\theta)\right)^{-1}
\nabla_{\theta}J(\theta)\right) = 0 \right. \right\}
\end{equation}
\end{thm}

\section{Experimental results}
\label{sec_simulation}

The simulations are performed using
a two-node network of {$M/G/1$} queues with feedback,
which is a similar setting as the one considered by Bhatnagar~\cite{Bhatnagar_2007_jour_TOMACS}.

\begin{table*}[ht]
\centering
\begin{tabular}{|c||c|c||c|c||c|c|}
\hline
& \multicolumn{2}{c||}{$\beta=0.01$}	& \multicolumn{2}{c||}{$\beta=0.05$}	& \multicolumn{2}{c|}{$\beta=0.25$}	\\
\cline{2-7}													
$q$	&	G$q$-SF2	&	N$q$-SF2	&	G$q$-SF2	&	N$q$-SF2	&	G$q$-SF2	&	N$q$-SF2\\
\hline													
\hline													
0.001	&	0.6680$\pm$0.0645	&	0.7875$\pm$0.1334	&	0.5621$\pm$0.0519	&	0.5772$\pm$0.0793	&	0.7531$\pm$0.0640	&	0.5191$\pm$0.0653\\
0.200	&	0.6598$\pm$0.0623	&	0.7577$\pm$0.0743	&	0.5355$\pm$0.0799	&	0.4527$\pm$0.0864	&	0.6984$\pm$0.1159	&	0.5011$\pm$0.0571\\
0.400	&	0.6736$\pm$0.0476	&	0.7026$\pm$0.0895	&	0.5477$\pm$0.0736	&	0.4169$\pm$0.0700	&	0.7140$\pm$0.0800	&	0.4630$\pm$0.0565\\
0.600	&	0.6202$\pm$0.0728	&	0.7083$\pm$0.0928	&	0.5475$\pm$0.0411	&	0.4418$\pm$0.0623	&	0.7178$\pm$0.0697	&	0.4578$\pm$0.0732\\
0.800	&	0.5909$\pm$0.0533	&	0.6796$\pm$0.0653	&	0.5605$\pm$0.0721	&	0.4256$\pm$0.0749	&	0.6427$\pm$0.0676	&	0.4475$\pm$0.0525\\
{Gaussian}&	0.6339$\pm$0.0658	&	0.6657$\pm$0.0816	&	0.5018$\pm$0.0647	&	0.4111$\pm$0.0534	&	0.6922$\pm$0.0670	&	0.4568$\pm$0.0635\\
1.020	&	0.6394$\pm$0.0738	&	0.6978$\pm$0.0732	&	0.4755$\pm$0.0701	&	0.4266$\pm$0.0685	&	0.7135$\pm$0.0763	&	0.4427$\pm$0.0518\\
1.040	&	0.6101$\pm$0.0663	&	0.6323$\pm$0.0768	&	0.4646$\pm$0.0405	&	0.3950$\pm$0.0807	&	0.6483$\pm$0.0514	&	0.4438$\pm$0.0602\\
1.060	&	0.6362$\pm$0.1036	&	0.6675$\pm$0.0894	&	0.4988$\pm$0.0796	&	0.3894$\pm$0.0520	&	0.7143$\pm$0.0755	&	0.4775$\pm$0.0556\\
1.080	&	0.5319$\pm$0.0745	&	0.6598$\pm$0.0787	&	0.5019$\pm$0.0353	&	0.4068$\pm$0.0503	&	0.6611$\pm$0.0866	&	0.4865$\pm$0.0680\\
{Cauchy}&	0.6217$\pm$0.0533	&	0.6455$\pm$0.0925	&	0.5359$\pm$0.0255	&	0.4573$\pm$0.0570	&	0.7040$\pm$0.0693	&	0.4861$\pm$0.0588\\
1.099	&	0.6440$\pm$0.0635	&	0.6577$\pm$0.0721	&	0.6550$\pm$0.1071	&	0.5722$\pm$0.0977	&	0.8658$\pm$0.1703	&	0.5873$\pm$0.0942\\
\hline													
\end{tabular}
\caption{$\Vert\theta(n)-\bar{\theta}\Vert$ for G$q$-SF2 and N$q$-SF2 for 
varying $q$ and $\beta$, when step-sizes are $a(n)=\frac{1}{(n+1)}$, $b(n) = \frac{1}{(n+1)^{0.85}}$
$c(n) = \frac{1}{(n+1)^{0.65}}$.}
\label{tab_compare}
\end{table*}

\begin{figure}[h]
\centering
\includegraphics[width=0.28\textwidth]{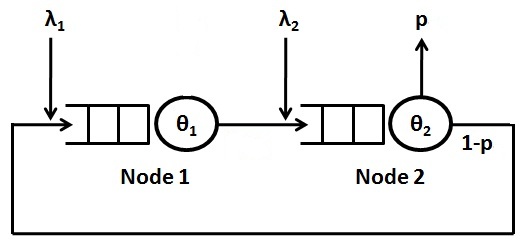}
\caption{Queuing Network.}
\label{fig_queue}
\end{figure}
The nodes in the network, shown in Fig.~\ref{fig_queue} are 
fed with independent Poisson external arrival processes with rates {$\lambda_1$} and {$\lambda_2$},
respectively. After service at the first node, a customer enters the second node.
When a customer departs from the second node, he either leaves the system with probability {$p = 0.4$} or 
re-enters the first node with the remaining probability. 
The service time processes at each node, {$\{S_n^i(\theta_i)\}_{n\geqslant1}, i = 1,2$} are given by
\begin{equation}
S_n^i(\theta_i) = U_i(n)\frac{\left(1+\Vert{\theta_i(n) - \bar{\theta}_i}\Vert^2\right)}{R_i} \;,
\end{equation}
where {$R_i$} are constants and {$U_i(n)$} are independent samples drawn from the uniform distribution on {$(0,1)$}. 
The service time of each node depends on the {$N_i$}-dimensional tunable parameter vector {$\theta_i$}, whose individual components
lie in the closed interval $[\alpha_{\min},\alpha_{\max}] = [0.1,0.6]$. 
{$\theta_i(n)$} represents the {$n^{th}$} update of the parameter vector at the {$i^{th}$} node, and {$\bar{\theta}_i$}
represents the target parameter vector corresponding to the {$i^{th}$} node. 
For the purpose of simulations, we consider $\lambda_1 = 0.2$, $\lambda_2 = 0.1$, $R_1 = 10$ and $R_2 = 20$.

The cost function, at any instant, is the total waiting time of all the customers in the system. 
In order to minimize the cost, we need to minimize {$S_n^i(\theta_i)$},
\textit{i.e.}, we require {$\theta_i(n)=\bar{\theta}_i$}, {$i=1,2$}. 
Let {$N=N_1+N_2$} and we consider $\theta,\bar{\theta} \in\mathbb{R}^N$ as
{$\theta = (\theta_1, \theta_2)^{T}$} and {$\bar\theta = (\bar\theta_1, \bar\theta_2)^{T}$}. 
Thus, $\bar{\theta}$ is the optimal value, and hence, we use
$\Vert{\theta(n)-\bar{\theta}}\Vert$ as a measure of performance of the algorithm.
The service time parameters at each node are assumed to be 10-dimensional vectors
($N_1=N_2=10$).
Thus, $N=20$ and {$C=[0.1,0.6]^{20}$}.
We fix each component of the target parameter vector, {$\bar{\theta}$}, at 0.3
and each component of the initial parameter, $\theta(0)$, at 0.6.
The simulations were performed using C on an Intel Pentium dual core machine 
with Linux operating system.

The analysis along the faster timescale for the Newton based algorithms shows that 
the gradient and Hessian updates run independently
and are not coupled between themselves, \textit{i.e.}, 
update of one does not influence the other,
and hence, their convergence to the smoothed gradient and Hessian, 
respectively, can be independently analyzed. 
This also provides a scope to update the gradient and Hessian along different timescales without affecting the
convergence of the algorithms. 
The step-size sequences for the parameter update and gradient estimation are 
chosen as $a(n) = \frac{1}{(n+1)}$ and $b(n) = \frac{1}{(n+1)^{0.85}}$,
respectively, while the one for Hessian estimation is considered as
$c(n) = \frac{1}{(n+1)^\gamma}$, $n\geqslant0$.
In order to satisfy Assumption~\ref{stepsize}, we require $\gamma\in(0.5,1)$,
but Bhatnagar~\cite{Bhatnagar_2007_jour_TOMACS} 
observed that better performance can be achieved in the N-SF2 algorithm if
Hessian is updated on a faster timescale. 
Even though as suggested by the convergence analysis, one does not
require three separate timescales, as two timescales are sufficient,
it is observed empirically that updating Hessian on a timescale
faster than both the parameter and the gradient updates can lead to better performance.

\begin{figure}[ht]
\centering
\includegraphics[width=0.37\textwidth]{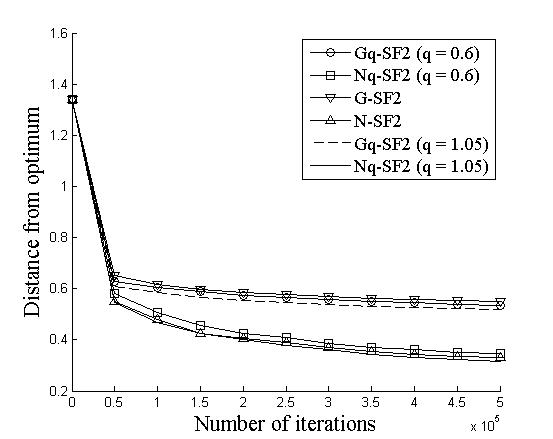}
\caption{Convergence behavior of Gaussian and $q$-Gaussian SF algorithms for {$q=0.6$}
and $1.05$.}
\label{plot}
\end{figure}

We compare the performance of the Jacobi variant of the N$q$-SF2 algorithm
with respect to the corresponding gradient based method (G$q$-SF2) for different 
values of $q$, $\beta$ and $\gamma$.
The other parameters are held fixed at {$M = 5000$}, {$L=100$} and {$\varepsilon=0.1$}. 
Thus, we perform a total of $2ML = 10^6$ simulations
to obtain $M=5000$ parameter updates. The following results are averaged over 
20 independent runs, each requiring about 5 seconds of clock time.
Fig.~\ref{plot} shows the convergence
behavior of the G$q$-SF2 and N$q$-SF2 for {$q=0.6$}, $q=1.05$ and Gaussian
with {$\beta=0.1$} and $\gamma=0.65$.

Table~\ref{tab_compare} compares the performance of
 G$q$-SF2 and N$q$-SF2, for different
values of $q$, in terms of the mean and variance of the final distance from the target vector. 
We note here that, in this case, $q$ varies in the range $(0,1.1)$
since $N=20$. The two special cases of Gaussian and Cauchy are retrieved
for $q\to1$ and $q = 1.095$, respectively. The table presents a comparison 
for three values of $\beta$, viz., $\beta = 0.01, 0.05$ and $0.25$, respectively (comparison for $\beta=0.1$ 
is given in the second and fourth columns of Table~\ref{tab_step}).
The step-size for the Hessian update is fixed with $\gamma=0.65$.
The results show that although for small $\beta$ ($\beta=0.01$), G$q$-SF2
works better than N$q$-SF2, but for higher $\beta$, N$q$-SF2 consistently outperforms
its gradient counterpart. In fact, it can be observed that the ratio of the distances
obtained using G$q$-SF2 and N$q$-SF2 increases with increasing $\beta$,
indicating that the relative performance of N$q$-SF2
in relation to G$q$-SF2 improves with more smoothing.
Other observations pertaining to the trends of performance with respect to $q$
and $\beta$ are similar to those for G$q$-SF2, discussed in~\cite{Ghoshdastidar_2013_arxiv}.

\begin{table*}[ht]
\centering
\begin{tabular}{|c||c|c|c||c|}
\hline
{$q$} & N$q$-SF2 $(\gamma=0.55)$ & N$q$-SF2 $(\gamma=0.65)$ & N$q$-SF2 $(\gamma=0.75)$ & G$q$-SF2	\\
\hline													
\hline													
0.001	&	0.4867$\pm$0.1056	&	0.4698$\pm$0.0627	&	\textbf{0.4335$\pm$0.0693}	&	0.5561$\pm$0.0832	\\
0.200	&	0.3847$\pm$0.0729	&	\textbf{0.3589$\pm$0.0591}	&	0.3698$\pm$0.0695	&	0.5555$\pm$0.0544	\\
0.400	&	\textbf{0.3328$\pm$0.0554}	&	0.3547$\pm$0.0548	&	0.3956$\pm$0.0698	&	0.5376$\pm$0.0569	\\
0.600	&	0.3422$\pm$0.0792	&	\textbf{0.3163$\pm$0.0582}	&	0.3254$\pm$0.0488	&	0.5472$\pm$0.0470	\\
0.800	&	\textbf{0.3127$\pm$0.0694}	&	0.3136$\pm$0.0699	&	0.3397$\pm$0.0536	&	0.5068$\pm$0.0548	\\
{Gaussian}	&	\textbf{0.3130$\pm$0.0539}	&	0.3560$\pm$0.0488	&	0.3383$\pm$0.0514	&	0.5354$\pm$0.0810	\\
1.020	&	\textbf{0.3160$\pm$0.0534}	&	0.3223$\pm$0.0397	&	0.3712$\pm$0.0653	&	0.5263$\pm$0.0511	\\
1.040	&	0.3203$\pm$0.0585	&	\textbf{0.3081$\pm$0.0552}	&	0.3315$\pm$0.0655	&	0.5103$\pm$0.0965	\\
1.060	&	\textbf{0.3130$\pm$0.0599}	&	0.3216$\pm$0.0566	&	0.3681$\pm$0.0540	&	0.4725$\pm$0.0599	\\
1.080	&	0.3722$\pm$0.0516	&	\textbf{0.3584$\pm$0.0633}	&	0.3782$\pm$0.0384	&	0.5165$\pm$0.0666	\\
{Cauchy}	&	0.4249$\pm$0.0615	&	\textbf{0.3997$\pm$0.0509}	&	0.4402$\pm$0.0657	&	0.5685$\pm$0.0798	\\
1.099	&	\textbf{0.5450$\pm$0.0683}	&	0.5594$\pm$0.0623	&	0.5677$\pm$0.0626	&	0.7253$\pm$0.0776	\\
\hline
\end{tabular}
\caption{$\Vert\theta(n)-\bar{\theta}\Vert$ for G$q$-SF2 and N$q$-SF2 for 
varying $q$ and varying step-size for Hessian update, $c(n) =\frac{1}{(n+1)^\gamma}$, with $\beta=0.1$ and other step-sizes 
maintained at $a(n)=\frac{1}{(n+1)}$ and $b(n) = \frac{1}{(n+1)^{0.85}}$.}
\label{tab_step}
\end{table*}

We also discuss about the effect of updating the Hessian estimate along different timescales.
Table~\ref{tab_step} shows the effect of {$\gamma$} on the N$q$-SF2 algorithm for varying {$q$}
(see earlier discussion),
while {$\beta$} is held fixed at {$0.1$}. It can be observed that, at this 
level of smoothing, N$q$-SF2 always performs better than G$q$-SF2. 
For each value of $q$, the best value of $\gamma$ is highlighted.
A faster update of the Hessian is seen to result in an improved performance.
Finally, it is interesting to note from both tables that the best results
are most often obtained for values of $q$ that do not correspond to 
either the Gaussian or the Cauchy perturbations,
thereby signifying the importance of generalization of the SF algorithms to include 
$q$-Gaussian perturbations, with a continuously-valued $q$ parameter.

\section{Conclusions}
\label{sec_conclusion}
We proposed a two-simulation SF algorithm with $q$-Gaussian perturbations
to perform Newton based optimization of a stochastic objective function.
In this process, we derived estimates for the Hessian of 
a two-sided smoothed functional using $q$-Gaussian distribution.
We also derived conditions for convergence of the algorithms,
and illustrated the performance of the algorithms through numerical simulations.

An interesting fact here is that though it is known 
that gradient of $q$-Gaussian SF always converges to the gradient of cost function 
as smoothing parameter $\beta\to0$, 
we observed that the same does not always hold for the Hessian.
In particular, we found that the Hessian in the case of uniform ($q\to-\infty$) smoothing
does not converge.
The issue lies in the attempt to derive an expression 
for the smoothed Hessian, $\nabla_\theta^2 S_{q,\beta}[J(\theta)]$,
in terms of the cost function. One can verify that this is not possible 
for the uniform case since the Hessian of the smoothed functional turns out to be
in the form of a finite difference of the gradient of the cost function.
Deriving a Hessian estimator for the case of uniform perturbations 
in terms of objective function remains an open problem.

As suggested in~\cite{Bhatnagar_2007_jour_TOMACS}, we may vary the 
smoothing parameter, $\beta$, at different update iterations. 
It would be useful to use more smoothing (larger $\beta$) at the initial
stages to proceed towards the global minimum, whereas at later stages of the algorithm
a smaller value of $\beta$ would provide better estimates for the gradient and Hessian.
From the analysis point of view, such a modification does not affect
the results, where we can easily 
replace $\beta$ by the corresponding sequence $\beta(n)$
as long as the sequence $\frac{b(n)}{\beta(n)}$ satisfies the conditions in Assumption~\ref{stepsize}
in place of the sequence $b(n)$. Further,
Theorem~\ref{thm_NqSF2} holds as long as we have $\sup_n \beta(n) \leqslant \beta_0$.

A similar modification may also be used for the values of $q$. It has been empirically observed,
both in Section~\ref{sec_simulation} and in~\cite{Ghoshdastidar_2013_arxiv},
that as $\beta$ decreases, larger values of $q$ tend to perform better. Hence,
one may start from a high value of $\beta$ and low value of $q$, and  
can decrease the former and increase the latter as the iterations proceed.
One may also incorporate the modification suggested in~\cite{Spall_2000_jour_AutoCtrl},
where steepest descent is employed for the initial parameter updates 
and Newton based search is employed for faster convergence of later recursions.
Such algorithmic modifications can be effectively used to improve the performance of the method 
and reduce the the computational burden of full Newton methods
without affecting the theoretical analysis.

\appendix
\section{Details of convergence analysis}
\label{app_convergence}

The convergence analysis is based on three key results. 
The convergence of the gradient and Hessian
recursions follow from results in~\cite{Borkar_2008_book_Cambridge}
(summarized in Lemma~\ref{borkar_cor_8}). While the main result
in our context is given in Lemma~\ref{W_bounded}, Lemma~\ref{B_n_convergence} 
proves a necessary condition required for the application of Lemma~\ref{borkar_cor_8}.

The main theorem for convergence of N$q$-SF2 is based on a result on projected
iterated schemes~\cite{Kushner_1978_book_Springer}. The conditions
in this result require Lemma~\ref{hess_xi_n_convergence} to hold.
While the result in~\cite{Kushner_1978_book_Springer} 
helps to eliminate the effect of the noise term,
Proposition~\ref{hess_SF2_convergence} shows that the error term is also small
and does not affect convergence. A more rigorous way to use
the consequences of Proposition~\ref{hess_SF2_convergence} would be via 
Hirsch's lemma~\cite{Hirsch_1989_jour_NeuNet}. We provide the intuitive 
arguments for this in the proof of Theorem~\ref{thm_NqSF2}.

Lastly, most of the proofs given below use a technical result regarding 
statistical properties of the 
multivariate $q$-Gaussian distribution~\cite[Proposition~4.1]{Ghoshdastidar_2013_arxiv}.
This result provides 
conditions for existence and an expression for the following expectation
\begin{equation}
\mathsf{E} \left[\frac{\left(\eta^{(1)}\right)^{b_1}\left(\eta^{(2)}\right)^{b_2}\ldots\left(\eta^{(N)}\right)^{b_N}}{
\left(\rho(\eta)\right)^b}\right]
\label{eq_gen_moment}
\end{equation}
for a standard $N$-dimensional $q$-Gaussian random variate $\eta = (\eta^{(1)},\ldots,\eta^{(N)})$ 
and non-negative integers $b$, $b_1, \ldots, b_N$, where $\rho$ is as defined in~\eqref{rho_defn}.
We skip the details of this result, but state few 
consequences in the following corollary that will be used in later discussions.
The claims below immediately follow from~\cite[Proposition~4.1]{Ghoshdastidar_2013_arxiv}.

\begin{cor}
 \label{cor_moments}
 The expectation in~\eqref{eq_gen_moment} exists and is finite  
 for  {$b<\big(1+\frac{1}{1-q}\big)$} when {$q<1$}, and the same holds for 
 {$b > \left(\frac{N}{2} - \frac{1}{q-1} +\sum_{i=1}^{N}\frac{b_i}{2}\right)$} if {$1<q<\left(1+\frac{2}{N}\right)$}.
 Further, in special cases, we have the following simplifications.
 Using notations similar to~\eqref{eq_gen_moment}, we have
 \begin{enumerate}
  \item 
  the term in~\eqref{eq_gen_moment} is zero whenever at least one of the $b_i$'s is odd,
  \item 
  $\mathsf{E} \displaystyle\left[ \frac{{\eta^{(i)}}^2}{\rho(\eta)} \right] = 
  \mathsf{E} \left[ \frac{1}{\rho(\eta)} \right] = \frac{N+2-Nq}{2}$,
  \item
  $\mathsf{E} \displaystyle\left[ \frac{{\eta^{(i)}}^2}{\rho(\eta)^2} \right] = \frac{N+2-Nq}{4q}$, and
  \item
  $\mathsf{E} \displaystyle\left[ \frac{{\eta^{(i)}}^4}{\rho(\eta)^2} \right] = 
  3\mathsf{E} \left[ \frac{{\eta^{(i)}}^2{\eta^{(j)}}^2}{\rho(\eta)^2} \right] = \frac{3(N+2-Nq)^2}{4q}$
 \end{enumerate}
 for all $i,j = 1,\ldots,N$, $i \neq j$.
 The latter two statements involve a term of $\rho(\eta)^2$, and hence,
 exist only for $q\in(0,1)\cup(1,1+\frac{2}{N})$, while the first two are defined over the entire range of $q$'s.
\end{cor}

\subsection*{Detailed proofs of results in Section~\ref{sec_convergence}}
The proofs of Lemmas~\ref{B_n_convergence} and~\ref{W_bounded} and Proposition~\ref{hess_SF2_convergence}
contain two parts related to the gradient and Hessian. We mostly prove the result
for the Hessian case. The corresponding proofs for gradient can be approached 
in a similar (in fact, simpler) manner, which are presented in~\cite{Ghoshdastidar_2013_arxiv}.

\begin{pf}(Proof of Lemma~\ref{B_n_convergence})

It is obvious that for all $p\geqslant1$, {$\mathsf{E}[B_{p}|\mathcal{G}_{p-1}] = 0$},
which implies {$(B_p,\mathcal{G}_p)_{p\in\mathbb{N}}$} is a martingale difference sequence. 
Using Jensen's inequality, we have
\begin{align}
\mathsf{E}&\left[\left.\Vert B_{p}\Vert^2\right|\mathcal{G}_{p-1}\right]
\leqslant \frac{8}{\beta^4(N+2-Nq)^2} \times
\nonumber \\
&\qquad
\mathsf{E}\left[\left.\Vert H(\tilde{\eta}(p))\Vert^2 {\big(h^2(Y_p)+h^2(Y'_p)\big)}\right| \mathcal{G}_{p-1}\right],
\label{Jensen_H1}
\end{align}
where we use {$\Vert.\Vert$} to denote the 2-norm for the matrices {$B_p$} and 
{$H(\tilde{\eta}(p))$} for {$p\in\mathbb{N}$}.
Denoting the Frobenius norm by {$\Vert.\Vert_F$}, we can use the definition of $H(.)$ to write
\begin{align}
\Vert H(\eta) \Vert_F^2 
&=
\frac{4q^2 \Vert\eta\Vert^4}{(N+2-Nq)^2\rho(\eta)^4} 
\nonumber \\
&- \frac{4q \Vert\eta\Vert^2}{(N+2-Nq)\rho(\eta)^3} + \frac{N}{\rho(\eta)^2}\;.
\label{H_norm}
\end{align}

\noindent
For {$q\in(-\infty,1)$}, we use Holder's inequality and the fact that {$\Vert H(\eta) \Vert \leqslant \Vert H(\eta) \Vert_F$} to claim
\begin{align*}
\mathsf{E}&\left[\left.\Vert B_{p}\Vert^2\right|\mathcal{G}_{p-1}\right]
\leqslant \frac{8}{\beta^4(N+2-Nq)^2} \times
\nonumber \\
&\qquad
\sup\limits_{\eta} \left(\Vert H(\eta) \Vert_F^2\right)
\mathsf{E}\left[\left.{h^2(Y_p)+h^2(Y'_p)}\right| \mathcal{G}_{p-1}\right],
\end{align*}
where, from~\eqref{H_norm}, we can argue that {$\sup_{\eta} \Vert H(\eta) \Vert_F^2$} is finite 
for any finite $q\in(-\infty,1)$ as
{$0\leqslant\Vert\eta\Vert^2 < \frac{N+2-Nq}{1-q}$} and {$\rho(\eta) \geqslant1$} for all {$\eta\in\Omega_q$}. Further
by the Lipschitz continuity of $h$ and Assumption~\ref{lyapunov}, we can claim 
{$\mathsf{E} \left[h^2(Y_p)|\mathcal{G}_{p-1}\right] <\infty$} and 
{$\mathsf{E} \left[h^2(Y'_p)|\mathcal{G}_{p-1}\right] <\infty$} a.s.
Thus, {$\mathsf{E}\left[\left.\Vert B_{p}\Vert^2\right|\mathcal{G}_{p-1}\right] < \infty$} a.s. 
for all {$p\in\mathbb{N}$}.
For {$q\in\big(1,1+\frac{2}{N}\big)$}, we note the second term in~\eqref{H_norm} is negative,
and hence, we may bound using only the first and third terms as
\begin{align}
&\mathsf{E}\left[\left.\Vert B_{p}\Vert^2\right|\mathcal{G}_{p-1}\right]
\nonumber \\
&\leqslant \frac{32q^2}{\beta^2(N+2-Nq)^4}
\mathsf{E}\left[\left.\frac{\Vert\eta\Vert^4}{\rho(\eta)^4} {\big(h^2(Y_p)+h^2(Y'_p)\big)}\right| \mathcal{G}_{p-1}\right]
\nonumber \\
&+ \frac{8N}{\beta^2(N+2-Nq)^2} 
\mathsf{E}\left[\left.\frac{\big(h^2(Y_p)+h^2(Y'_p)\big)}{\rho(\eta)^2} \right| \mathcal{G}_{p-1}\right].
\label{Jensen_H2}
\end{align}
Applying the Cauchy-Schwartz inequality 
and the fact that $(a+b)^2 \leqslant 2(a^2+b^2)$, for any $a,b\in\mathbb{R}$,
on the first term in~\eqref{Jensen_H2}, we obtain
\begin{align*}
&\mathsf{E}\left[\left.\frac{\Vert\eta\Vert^4}{\rho(\eta)^4} {\big(h^2(Y_p)+h^2(Y'_p)\big)}\right| \mathcal{G}_{p-1}\right]
\\&\leqslant 
\sqrt{2}\mathsf{E}\left[\frac{\Vert\eta\Vert^8}{\rho(\eta)^8}\right]^{1/2}
\mathsf{E}\left[\left.{\big(h^4(Y_p)+h^4(Y'_p)\big)}\right| \mathcal{G}_{p-1}\right]^{1/2}.
\end{align*}
The second expectation  is finite a.s. from earlier discussion. 
We expand {$\Vert\eta\Vert^8$} in the first expectation
and use Corollary~\ref{cor_moments}
to claim the existence and finiteness of the expectation for {$q\in\big(1,1+\frac{2}{N}\big)$}.
Similar arguments are applicable for the second term in~\eqref{Jensen_H2}
and the claim follows. 
\end{pf}

\begin{pf}(Proof of Lemma~\ref{W_bounded})

Since Lemma~\ref{B_n_convergence} holds,
one can verify that iterations~\eqref{fastode_Z} and~\eqref{fastode_W}
satisfy the necessary conditions
required to apply Lemma~\ref{borkar_cor_8},
where the invariant measure, $\nu$, of the process $\{(Y_p,Y'_p)\}_p$ is
the product measure of $\nu_{(\theta+\beta\eta)}$ and $\nu_{(\theta-\beta\eta)}$, 
the invariant measures of the processes $\{Y_p\}_p$ and $\{Y'_p\}_p$, respectively.
The claim follows from an application of the aforementioned result.
\end{pf}

\begin{pf}(Proof of Proposition~\ref{hess_SF2_convergence})

For small {$\beta>0$}, we use Taylor's expansion of 
{$J(\theta+\beta\eta)$} and {$J(\theta-\beta\eta)$} around {$\theta\in C$} to write
\begin{align*}
J(\theta+\beta\eta) &+ J(\theta-\beta\eta) \nonumber
\\&= 2J(\theta) + \beta^2\eta^{T}\nabla_{\theta}^2 J(\theta)\eta + o(\beta^3).
\\
J(\theta+\beta\eta) &- J(\theta-\beta\eta)
= 2\beta\eta^{T}\nabla_{\theta}J(\theta) + o(\beta^2).
\end{align*}
Thus the gradient of the two-sided SF~\eqref{grad_SF2_formula} becomes
\begin{align}
&\nabla_{\theta}S_{q,\beta}[J(\theta)] 
\nonumber \\
&=\frac{1}{(N+2-Nq)}
\mathsf{E}_{G_q(\eta)}\left[\frac{2}{\rho{(\eta)}}\eta\eta^T\right]\nabla_{\theta}J(\theta) + o(\beta),
\end{align}
and the two-sided smoothed Hessian~\eqref{hess_SF2_formula} is
\begin{align}
\nabla_{\theta}^2 S_{q,\beta}&[J(\theta)]
= \frac{1}{\beta^2 (N+2-Nq)} \bigg(
2J(\theta)\mathsf{E}\left[H(\eta)|\theta\right] 
\nonumber\\
&+\beta^2\mathsf{E}\left[H(\eta)\eta^{T}\nabla_{\theta}^2J(\theta)\eta\,|\theta\right] +
o(\beta^3) \bigg).
\label{hess_expand}
\end{align}
Let us consider each of the terms in~\eqref{hess_expand}. 
Corollary~\ref{cor_moments} ensures that 
the product moments are zero whenever the product is odd.
Hence, for all {$i,j=1,\ldots,N$}, {$i\neq j$}, 
{$\mathsf{E}\left[H(\eta)_{i,j}\right] = 0$}. 
Thus, the off-diagonal terms are zero, whereas the diagonal elements
are of the form
\begin{equation}
\mathsf{E}\left[H(\eta)_{i,i}\right]
= \frac{2q}{(N+2-Nq)}\mathsf{E}\left[\frac{\left(\eta^{(i)}\right)^2}{\rho{(\eta)}^2}\right]
- \mathsf{E}\left[\frac{1}{\rho{(\eta)}}\right],
\label{hess_expand_term1}
\end{equation}
for all {$i=1,2,\ldots,N$}.
Corollary~\ref{cor_moments}
shows that the expectations in~\eqref{hess_expand_term1} exist
for $q\in(0,1)\cup(1,1+\frac{2}{N})$. 
One can note that the squared term in the denominator imposes the condition $q>0$.
Substituting the corresponding expressions
in~\eqref{hess_expand_term1}, we get {$\mathsf{E}\left[H(\eta)_{i,i}\right] = 0$}.
Thus, the first term in~\eqref{hess_expand} is zero. 
Now, we consider the second term. For {$i\neq j$}, 
\begin{align*}
&\mathsf{E}\left[H(\eta)_{i,j} \left(\eta^{T}\nabla_{\theta}^2J(\theta)\eta\right)\,|\theta\right] 
\\&= \frac{2q}{(N+2-Nq)} \sum_{k,l=1}^{N}
\left[\nabla_{\theta}^2 J(\theta)\right]_{k,l} 
\mathsf{E}\left[\frac{\eta^{(i)}\eta^{(j)}\eta^{(k)}\eta^{(l)}}{\rho{(\eta)}^2}\right],
\end{align*}
which is zero unless {$i=k, j=l$} or {$i=l,j=k$}. So using the fact that {$\nabla_{\theta}^2 J(\theta)$} 
is symmetric, \textit{i.e.}, {$\left[\nabla_{\theta}^2 J(\theta)\right]_{k,l} = \left[\nabla_{\theta}^2 J(\theta)\right]_{l,k}$},
we can write
\begin{align}
&\mathsf{E}\left[H(\eta)_{i,j} \left(\eta^{T}\nabla_{\theta}^2J(\theta)\eta\right)\,|\theta\right] 
\nonumber \\
&= \frac{4q}{(N+2-Nq)} \left[\nabla_{\theta}^2 J(\theta)\right]_{i,j}
\mathsf{E}\left[\frac{\left(\eta^{(i)}\right)^2\left(\eta^{(j)}\right)^2}{\rho{(\eta)}^2}\right].
\label{hess_expand_term3_1}
\end{align}
Referring to Corollary~\ref{cor_moments}, 
we obtain
\begin{displaymath}
\mathsf{E}\left[H(\eta)_{i,j} \left(\eta^{T}\nabla_{\theta}^2J(\theta)\eta\right)\,|\theta\right]
=  (N+2-Nq)\left[\nabla_{\theta}^2 J(\theta)\right]_{i,j}
\end{displaymath}
for {$i\neq j$}.
Now for {$i=j$}, we use the definition of $H$~\eqref{H_defn} to write
\begin{align*}
&\mathsf{E}\left[H(\eta)_{i,i} \left(\eta^{T}\nabla_{\theta}^2J(\theta)\eta\right)\,|\theta\right]
\\&= \frac{2q}{(N+2-Nq)} \sum_{k,l=1}^{N}
\left[\nabla_{\theta}^2 J(\theta)\right]_{k,l} 
\mathsf{E}\left[\frac{\left(\eta^{(i)}\right)^2\eta^{(k)}\eta^{(l)}}{\rho{(\eta)}^2}\right]
\\&\qquad\qquad\qquad
- \sum_{k,l=1}^{N} \left[\nabla_{\theta}^2 J(\theta)\right]_{k,l} 
\mathsf{E}\left[\frac{\eta^{(k)}\eta^{(l)}}{\rho{(\eta)}}\right].
\end{align*}
Since the above expectations are zero for {$k\neq l$}, we have
\begin{align}
\mathsf{E}&\left[H(\eta)_{i,i} \left(\eta^{T}\nabla_{\theta}^2J(\theta)\eta\right)\,|\theta\right]
= \frac{2q \left[\nabla_{\theta}^2 J(\theta)\right]_{i,i}}{(N+2-Nq)} 
\mathsf{E}\left[\frac{\left(\eta^{(i)}\right)^4}{\rho{(\eta)}^2}\right] 
\nonumber \\
&+\frac{2q }{(N+2-Nq)}  \sum_{k\neq i}
\left[\nabla_{\theta}^2 J(\theta)\right]_{k,k} 
\mathsf{E}\left[\frac{\left(\eta^{(i)}\right)^2\left(\eta^{(k)}\right)^2}{\rho{(\eta)}^2}\right]  
\nonumber \\
&- \sum_{k=1}^{N} \left[\nabla_{\theta}^2 J(\theta)\right]_{k,k} 
\mathsf{E}\left[\frac{\left(\eta^{(k)}\right)^2}{\rho{(\eta)}}\right].
\label{hess_expand_term3_2}
\end{align}
We again refer to Corollary~\ref{cor_moments}
to compute each term in~\eqref{hess_expand_term3_2},
and then perform simple algebraic manipulations to derive
\begin{align*}
\mathsf{E}\left[H(\eta)_{i,i} \left(\eta^{T}\nabla_{\theta}^2J(\theta)\eta\right)|\theta\right]
= (N+2-Nq) [\nabla_{\theta}^2 J(\theta)]_{i,i}.
\end{align*}
By substituting all the above expressions in~\eqref{hess_expand},
we have that the difference between $\nabla_{\theta}^2 S_{q,\beta}[J(\theta)]$ and $\nabla_{\theta}^2 J(\theta)$
is $o(\beta)$, which implies that the Euclidean distance 
between $\nabla_{\theta}^2 S_{q,\beta}[J(\theta)]$ and $\nabla_{\theta}^2 J(\theta)$ is $o(\beta)$.
A similar result can be shown for the gradient as well. 

For the second part of the claim,
we write {$\Delta(\theta)=$} 
\begin{align*}
&\left(\mathcal{P}_{pd}\left(\nabla_{\theta}^2 J(\theta)\right)^{-1}
- \mathcal{P}_{pd}\left(\nabla_{\theta}^2 S_{q,\beta}[J(\theta)]\right)^{-1}\right)
\nabla_{\theta} J(\theta) 
\\&+ \mathcal{P}_{pd}\left(\nabla_{\theta}^2 S_{q,\beta}[J(\theta)]\right)^{-1}
\left(\nabla_{\theta} J(\theta) - \nabla_{\theta} S_{q,\beta}[J(\theta)]\right) ,
\end{align*}
which implies that
\begin{align}
&\Big\Vert\Delta (\theta)\Big\Vert \leqslant 
\nonumber \\
&\left\Vert \mathcal{P}_{pd}\left(\nabla_{\theta}^2 J(\theta)\right)^{-1}
- \mathcal{P}_{pd}\left(\nabla_{\theta}^2 S_{q,\beta}[J(\theta)]\right)^{-1} \right\Vert
\Big\Vert \nabla_{\theta} J(\theta) \Big\Vert \nonumber
\\&+
\left\Vert \mathcal{P}_{pd}\left(\nabla_{\theta}^2 S_{q,\beta}[J(\theta)]\right)^{-1} \right\Vert
\Big\Vert \nabla_{\theta} J(\theta) - \nabla_{\theta} S_{q,\beta}[J(\theta)] \Big\Vert \;.
\label{error_bound}
\end{align}
Since {$\nabla_{\theta} J(\theta)$} is continuously differentiable on the compact set $C$,
{$\sup_{\theta\in C}\Vert\nabla_{\theta} J(\theta)\Vert <\infty$}. 
Also, since {$\mathcal{P}_{pd}(A)$} is a positive definite
matrix for any $N\times N$ matrix $A$, its inverse always exists, \textit{i.e.},  
{$\Vert(\mathcal{P}_{pd}(A))^{-1}\Vert <\infty$} considering any matrix norm. Thus, in order to justify the claim, 
we need to show that other terms are {$o(\beta)$}.
From the first part of the claim, we have
{$\Vert \nabla_{\theta} J(\theta) - \nabla_{\theta} S_{q,\beta}[J(\theta)] \Vert=o(\beta)$},
and we can write
\begin{align*}
&\left\Vert \mathcal{P}_{pd}\left(\nabla_{\theta}^2 J(\theta)\right)^{-1}
- \mathcal{P}_{pd}\left(\nabla_{\theta}^2 S_{q,\beta}[J(\theta)]\right)^{-1} \right\Vert
\\&= \Big\Vert 
\mathcal{P}_{pd}\left(\nabla_{\theta}^2 J(\theta)\right)^{-1}
\mathcal{P}_{pd}\left(\nabla_{\theta}^2 S_{q,\beta}[J(\theta)]\right)^{-1}
\times \\&\qquad\qquad\qquad
\Big(\mathcal{P}_{pd}\left(\nabla_{\theta}^2 S_{q,\beta}[J(\theta)]\right)
-\mathcal{P}_{pd}\left(\nabla_{\theta}^2 J(\theta)\right)
\Big)\Big\Vert
\\&\leqslant
\left\Vert \mathcal{P}_{pd}\left(\nabla_{\theta}^2 J(\theta)\right)^{-1} \right\Vert
\left\Vert \mathcal{P}_{pd}\left(\nabla_{\theta}^2 S_{q,\beta}[J(\theta)]\right)^{-1} \right\Vert 
\times \\&\qquad\qquad\qquad
\Big\Vert \mathcal{P}_{pd}\left(\nabla_{\theta}^2 S_{q,\beta}[J(\theta)]\right)
-\mathcal{P}_{pd}\left(\nabla_{\theta}^2 J(\theta)\right) \Big\Vert \;.
\end{align*}
We note that for any matrix $A$, the eigenvalues of {$\mathcal{P}_{pd}(A)$} are lower bounded by {$\varepsilon>0$}. 
Hence, the first two terms, which are upper bounded by the maximum eigenvalues of the inverse of the projected matrices,
can at most be {$\frac{1}{\varepsilon}$}.  
Also we have shown
the third term is $o(\beta)$. The claim follows.
\end{pf}

\begin{pf}(Proof of Lemma~\ref{hess_xi_n_convergence})
 
As {$\theta(k)$} is {$\mathcal{F}_{k}$}-measurable, while {$\eta(k)$} is independent of {$\mathcal{F}_{k}$}
for all {$k\geqslant0$}, we can conclude that {$\mathsf{E}[\xi_{k}|\mathcal{F}_{k}] = 0$}.
Thus 
{$(M_k,\mathcal{F}_k)_{k\geqslant0}$} is a
martingale sequence. 
Now note as in  Lemma~\ref{B_n_convergence} that
\begin{align*}
&\mathsf{E}\big[\left.\left\Vert\xi_{k}\right\Vert^2\right|\mathcal{F}_{k}\big]
\\&\leqslant 4  \mathsf{E}\Bigg[
\Bigg\Vert\mathcal{P}_{pd}
\left(\frac{H(\eta(k))\bar{J}_k}{\beta^2 (N+2-Nq)}\right)^{-1} \times
\\& \qquad\qquad\qquad\qquad\qquad\qquad
\frac{\eta(k) \bar{J}'_k}{\beta\rho{(\eta(k))}(N+2-Nq)}\Bigg\Vert^2\Bigg]
\\&\leqslant \frac{4}{\beta^2 (N+2-Nq)^2} \times
\\&
\mathsf{E}\left[
\left\Vert\mathcal{P}_{pd}
\left(\frac{H(\eta(k))\bar{J}_k}{\beta^2 (N+2-Nq)}\right)^{-1}\right\Vert^2
\left\Vert\frac{\eta(k)}{\rho{(\eta(k))}}\right\Vert^2 (\bar{J}'_k)^2\right]
\\&\leqslant \frac{8}{\varepsilon^2\beta^2 (N+2-Nq)^2} \times
\\&
\mathsf{E}\left[
\left\Vert\frac{\eta(k)}{\rho{(\eta(k))}}\right\Vert^2
\left(J(\theta(k)+\beta\eta(k))^2+J(\theta(k)-\beta\eta(k))^2\right)\right]
\end{align*}
since the first term  in the expectation is square of 
the maximum eigenvalue of the inverse of the projected Hessian matrix,
which can be bounded above by {$\frac{1}{\varepsilon^2}$}.   
Using an argument similar to Lemma~\ref{B_n_convergence}, we can show
$\mathsf{E}[\Vert\xi_{k}\Vert^2|\mathcal{F}_{k}]$
has bounded variance. Using the fact that {$\sum_n a(n)^2 < \infty$},
we can write
\begin{align*}
\sum_{n=0}^{\infty}\mathsf{E}\left[\Vert{M_{n+1}-M_n}\Vert^2\right]
&\leqslant \sum_{n=0}^{\infty} a(n)^2 \sup_n \mathsf{E}\left[\Vert\xi_{n}\Vert^2\right]
\end{align*}
is finite a.s. From here, 
the claim follows from the martingale convergence theorem~\cite[page 111]{Williams_1991_book_Cambridge}.
\end{pf}

\begin{pf}(Proof of Theorem~\ref{thm_NqSF2})
Before proving the result, we briefly discuss the roles of Assumptions~\ref{differentiable}--\ref{projection}
in the proof. Assumption~\ref{differentiable} is essential for defining the smoothed 
gradient~\eqref{grad_SF2_formula} and Hessian~\eqref{hess_SF2_formula}, while Assumption~\ref{ergodic}
helps us to define the long run cost~\eqref{Jdefn}.
The  existence of the stochastic Lyapunov function (Assumption~\ref{lyapunov}) is used for the 
application of Theorem~\ref{borkar_cor_8}. The assumption on the step-sizes (Assumption~\ref{stepsize})
is required to ensure the second condition in Theorem~\ref{borkar_cor_8} and last two conditions
of Lemma~\ref{kushner_thm_5.3.1}. Here, we note that the square summability of the sequence $(a(n))_n$
makes it easy to satisfy the a.s. convergence of the martingale sequence (c.f. Lemma~\ref{hess_xi_n_convergence}).
Finally, due to Assumption~\ref{projection}, one can prove 
the claims in Proposition~\ref{hess_SF2_convergence} and Lemma~\ref{hess_xi_n_convergence}, and hence, we
can apply Kushner and Clark's result,
see~\cite[Theorem~5.3.1]{Kushner_1978_book_Springer} 
(alternatively~\cite[Lemma~4.6]{Ghoshdastidar_2013_arxiv}),
to claim that 
the update in~\eqref{hess_slowstep} converges to the stable fixed points of the ODE
\begin{equation}
\label{slowode_error}
\dot{\theta}(t) = \tilde{\mathcal{P}}_C \Big(-\nabla_{\theta}J(\theta(t)) + \Delta\big(\theta(t)\big)\Big)\;.
\end{equation}
Note that the other condition in Lemma~\ref{kushner_thm_5.3.1} follows from the definition of $C$ and $\tilde{P}_C$,
and the continuous differentiability of $J$.

Now, starting from the same initial condition, 
if {$\Delta(\theta(t))\to0$}, then
the trajectory of~\eqref{slowode_error} 
tracks the trajectory of the following ODE:
\begin{equation}
\dot{\theta}(t) = \tilde{\mathcal{P}}_C \Big(-\nabla_{\theta}J(\theta(t))\Big)
\end{equation}
uniformly over compacts.
From Proposition~\ref{hess_SF2_convergence}, 
we have {$\left\Vert\Delta\big(\theta(n)\big)\right\Vert = o(\beta)$} for all $n$,
which leads to the claim.
\end{pf}


\end{document}